\documentclass[12pt]{article}	

\usepackage{amscd,amsmath,amsxtra,amssymb,amsthm,latexsym,amsfonts,color,graphics}
\usepackage{setspace}

\newtheorem{Thm}{Theorem}[section]
\newtheorem{Lemma}[Thm]{Lemma}
\newtheorem{Prop}[Thm]{Proposition}
\newtheorem{Cor}[Thm]{Corollary}

\newtheorem{Statement}[Thm]{Statement}

\numberwithin{equation}{section}

\newcommand{\R}{{\mathbb R}}
\renewcommand{\O}{{\cal O}}

\newcommand{\C}{{\mathbb C}}

\newcommand{\N}{{\mathbb N}}

\newcommand{\re}{{\rm Re \ }}

\renewcommand{\Re}{{\rm Re \ }}

\newcommand{\MLS}{{\text{\large \rm ML$\Sigma$}}}

\newcommand{\Log}{{\rm Log \ }}

\newcommand{\var}{{\rm var \ }}

\author{Alexander GETMANENKO \\ 
\footnotesize Department of Mathematics, Northwestern University, Evanston IL, U.S.A.;\\  
\footnotesize
Max Planck Institute for Mathematics in the Sciences, Leipzig, Germany; \\
\footnotesize Institute for the Physics and Mathematics of the Universe, \\
\footnotesize The University of Tokyo, 5-1-5 Kashiwanoha, Kashiwa, 277-8568, Japan\\
\footnotesize e-mail: {\tt Alexander.Getmanenko@ipmu.jp}; \ \ fax: +81-4-7136-4941 }

\title{On eigenfunctions corresponding to a small resurgent eigenvalue.}

\begin{document}

\maketitle		

\begin{abstract}  The paper is devoted to some foundational questions in resurgent analysis. As a main technical result, it is shown that under appropriate conditions the infinite sum of endlessly continuable majors commutes with the Laplace transform. A similar statement is proven for compatibility of a convolution and of an infinite sum of majors. We generalize the results of Candelpergher-Nosmas-Pham and prove a theorem about substitution of a small (extended) resurgent function into a holomorphic parameter of another resurgent function. Finally, we discuss an application of these results to the question of resurgence of eigenfunctions of a one-dimensional Schr\"odinger operator corresponding to a small resurgent eigenvalue. 
\end{abstract}

\newpage

\section{Introduction}

A resurgent function $f(h)$ is, roughly, a function admitting a nice enough hyperasymptotic expansion 
$$\sum_k e^{-c_k/h}(a_{k,0}+a_{k,1}h+a_{k,2}h^2+...), \ \  \ \  h\to 0,$$
and resurgent analysis is a way of working with such expansions by means of representing $f(h)$ as a Laplace integral of a ramified analytic function in the complex plane.  

Applications of the theory of resurgent functions to quantum-mechanical problems  (~\cite{V83},~\cite{DDP97},~\cite{DP99}) have demonstrated its elegance and computational power, which brings forward the question of a fully rigorous mathematical justification of these methods.

This article is devoted to some foundational questions of resurgent analysis as applied to the Schr\"odinger equation in one dimension. Suppose, as we explicitly mention in ~\cite{G} or as it implicitly happens, say, in ~\cite{DDP97}, that using methods of  ~\cite{ShSt} one can construct a resurgent solution of the Schr\"odinger equation in one dimension
$$ [-h^2\partial_q^2 + V(q,h)]\psi(q,h) = 0 $$
for the potential $V(q,h)$ analytic with respect to $q$ and polynomial in $h$. This situation includes, in particular, the eigenvalue problem 
\begin{equation}
[-h^2\partial_q^2 + V(q,h)]\psi(E_1,q,h) = hE_1\psi(E_1,q,h),  \label{eq1}
\end{equation}
where $\psi$ is required or expected to holomorphically depend on $q$ and $E_1$ and to be a resurgent function of $h$ for any fixed $q$ and $E_1$. (We say ``expected to be" because a complete justification of these facts is still an open problem.)
Then, imposing appropriate boundary conditions on $\psi$ and using methods of asymptotic matching (see e.g. ~\cite{DP99,DDP97}) one can write the quantization conditions and obtain by solving them some energy level $E_1(h)$ which is expected to be a resurgent function of $h$ which is almost never a polynomial; thus, substituting $E_1(h)$ into \eqref{eq1} brings us outside of the class of equations polynomial in $h$. In this paper we show that if $E_1$ is a small resurgent function (see below for a precise definition), then one can substitute $E_1$ into a resurgent solution $\psi(E_1,q,h)$ of \eqref{eq1} and obtain another resurgent function of $h$ that would solve the corresponding differential equation (see section \ref{ApplicnToExistence}).

In order to work out the detail of substituting a small resurgent function of $h$ into a holomorphic parameter of another resurgent function, the following more technical issues are addressed.

A resurgent function is given as a Laplace integral of a ramified analytic function (called  ``major") along some infinite contour. In order for the Laplace integral to converge, one needs to correct that major by an entire function so that it stays bounded along the infinite branches of the contour; the details of this procedure are discussed in ~\cite{CNP}. If we are to deal with an infinite sum of majors, we need to choose the corresponding corrections by entire functions in a coherent way, in the precise sense of section \ref{CaseOfConvergSeriesOfMajors}. Then one can show that the Laplace transform is compatible, under certain conditions, with an infinite sum of majors. 

A similar statement is shown for a convolution and for the reconstruction homomorphism (the reconstruction of a major from its decomposition into microfunctions). 

We also consider the question of substituting of a small resurgent function into a holomorphic function and generalize the argument of ~\cite{CNP}  to the case of what they would call an extended resurgent function. This case is reduced to the case considered in ~\cite{CNP} by means of interchanging the infinite sums with convolutions and Laplace transforms as above.
The ~\cite{CNP}'s proof will also be recalled in a slightly generalized form that allows to accomodate the case of substituting of  $k$ small resurgent functions into a holomorphic functions of $k$ variables.

Although the results of this article are perfectly natural, they have not been, to our knowledge, properly documented in the literature and we felt it is important to write up their detailed proofs.


\section{Preliminaries from Resurgent Analysis}

We need to combine the setup of ~\cite{CNP} (the resurgence is with respect to the semiclassical parameter $h$ rather than the coordinate $q$) and mathematical clarity of ~\cite{ShSt} and therefore have to mix their terminology and to somewhat change their notation. 

\subsection{Laplace transform and its inverse. Definition of a resurgent function.}

Morally speaking, we will be studying analytic functions $\varphi(h)$ admitting asymptotic expansions $e^{-c/h}(a_0+a_1 h + a_2 h^2+...)$ for $h\to 0$ and $\arg h$ constrained to lie in an arc $A$ of the circle of directions on the complex plane; respectively, the inverse asymptotic parameter $x=1/h$ will tend to infinity and $\arg x$ will belong to the complex conjugate arc $A^*$.  Such functions can be represented as Laplace transforms of functions ${\bf \Phi}(\xi)$ where the complex variable $\xi$ is Laplace-dual to $x=1/h$ and ${\bf \Phi}(\xi)$ is analytic in $\xi$ for $|\xi|$ large enough and $\arg \xi \in {\check A}$, the copolar arc to $A$. The concept of a resurgent function will be defined by imposing conditions on analytic behavior of ${\bf \Phi}$.

\subsubsection{"Strict" Laplace isomorphism}

For details see ~\cite[Pr\'e I.2]{CNP}.

Let $A$ be a small (i.e. of aperture $<\pi$) arc in the circle of directions $S^1$. Denote by $\check A$ its copolar arc, $\check A = \cup_{\alpha\in A} {\check \alpha}$, where $\check \alpha$ is the open arc of length $\pi$ consisting of directions forming an obtuse angle with $\alpha$. 

Denote by ${\cal O}^\infty(A)$ the space of sectorial germs at infinity in direction $A$ of holomorphic functions and by ${\cal E}(A)$ the subspace of those that are of exponential type in the direction $A$, i.e. bounded by $e^{K|t|}$ as the complex argument $t$ goes to infinity in the direction $A$.
 
We want to construct the isomorphism
$$ {\cal L} \ : \ {\cal E}({\check A})/{\cal O}(\C)^{exp} \ \longleftrightarrow \ {\cal E}(A^*) \ : \ {\bar{\cal L}}$$ 
where ${\cal O}(\C)^{exp}$ denotes the space of functions of exponential type in all directions.

{\bf Construction of ${\cal L}$}. Let ${\bf \Phi}$ be a function holomorphic in a sectorial neighborhood $\Omega$ of infinity in direction $\check A$. For any small arc $A'\subset\subset A$ we can choose a point $\xi_0\in\C$ such that $\Omega$ contains a sector $\xi_0 {\check A}'$ with the vertex $\xi_0$ opening in the direction ${\check A}'$. Define the Laplace transform 
$$ \Phi_\gamma (x) := \int_\gamma e^{-x\xi} {\bf \Phi}(\xi) d\xi $$
with $\gamma=-\partial (\xi_0 {\check A}')$. Then $\Phi_\gamma$ is holomorphic of exponential type in a sectorial neighborhood of infinity in direction $A^*$. Cauchy theorem shows that if ${\bf \Phi}$ is entire of exponential type, then $\Phi_\gamma = 0$, which allows us to define
$$ {\cal L}({\bf \Phi} \ \rm{mod \ } {{\cal O}(\C)^{exp}} ) \ = \ \Phi_\gamma . $$

The construction of ${\bar {\cal L}}={\cal L}^{-1}$ will not be used in this paper.

\subsubsection{"Large" Laplace isomorphism}

The Laplace transform ${\cal L}$ defined in the previous subsection can be applied only to a function ${\bf \Phi}(\xi)$ satisfying a growth condition at infinity. At a price of changing the target space of ${\cal L}$ this restriction can be removed as follows.

For details see ~\cite[Pr\'e I.3]{CNP}.

Let ${\check A}=(\alpha_0,\alpha_1)$ be the copolar of a small arc, where $\alpha_0,\alpha_1\in S^1$ are two directions in the complex plane, and let $\gamma: \R\to \C$ be an endless continuous path. We will say that $\gamma$ is \underline{adapted} to ${\check A}$ if    $\lim_{t\to -\infty} \gamma(t)/|\gamma(t)|\to \alpha_0$, $\lim_{t\to \infty} \gamma(t)/|\gamma(t)|\to \alpha_1$, and if the length of the part of $\gamma$ contained in a ring $\{ z: R \le |z|\le R+1\}$ is bounded by a constant independent of $R$.

For a small arc $A$ there are two mutually inverse isomorphisms
$$ {\cal L} \ : \ {\cal O}^\infty({\check A})/{\cal O}(\C) \ \longleftrightarrow \ {\cal E}(A^*) / {\cal E}^{-\infty}(A^*) \ : \ {\bar{\cal L}},$$ 
where ${\cal E}^{-\infty}(A^*)$ is the set of sectorial germs at infinity that decay faster than any function of exponential type (cf. ~\cite[Pr\'e I.0]{CNP}).

{\bf Construction of ${\cal L}$}. Let ${\bf \Psi}$ be holomorphic in $\Omega$, a sectorial neighborhood of infinity of direction $\check A$. Let $\gamma$ be a path adapted to ${\check A}$ in $\Omega$. As we will see in lemma \ref{Lemma21} below, there is a function ${\bf \Phi}$ bounded on $\gamma$ such that ${\bf \Phi}-{\bf \Psi} \in {\cal O}(\C)$;   define
\begin{equation} {\cal L}({\bf \Psi} \mod {\cal O}(\C))  \ := \  \int_\gamma e^{-x\xi} {\bf \Phi}(\xi)d\xi \ \mod {\cal E}^{-\infty}(A^*). \label{defLfla} \end{equation}

{\bf Definition.} Any of the functions ${\bf \Psi}(\xi)$ satisfying ${\cal L}({\bf \Psi} {\rm\ mod \ } {\cal O}(\C)) = \psi(x) \mod {\cal E}^{-\infty}(A^*)$ is called a \underline{major} of the function $\psi(x)$.

An equivalence class of functions defined on a subset of $\C$ modulo adding an entire function is also called an \underline{integrality class}.

\subsubsection{Resurgent functions.}

Resurgent functions are usually understood to be functions of a large parameter $x$. For brevity we will speak of resurgent functions of $h$ to mean resurgent functions of $1/h$.  

{\bf Definition.} A germ $f(\xi)\in \O_{\xi_0}$  is \underline{endlessly continuable} if for any $L>0$ there is a finite set $\Omega_L\subset \C$ such that $f(\xi)$ has an analytic continuation along any path of length $<L$ avoiding $\Omega_L$.

{\bf Definition.} (cf. ~\cite[p.122]{ShSt})  Let $A\subset S^1\subset \C$ be a small arc and let $A^*$ be obtained from $A$ by complex conjugation. A \underline{resurgent function} $f(x)$ (of the variable $x\to\infty$) in direction $A^*$ is an element of ${\cal E}(A^*)/{\cal E}^{-\infty}$ such that any (hence all) of its majors is endlessly continuable.

{\bf Remark.} ~\cite{CNP} calls the same kind of object an "extended resurgent function". 

Let (cf. ~\cite[R\'es I]{CNP})  ${\pmb{\cal R}(A)}$ denote the set of endlessly continuable sectorial germs of analytic functions $\Phi(\xi)$ defined in a neighborhood of infinity in the direction $\check A$. Then a resurgent function of $x$ in the direction $A^*$ has a major in ${\pmb{\cal R}}(A)$. 

When we mean a resurgent function $g(h)$ of a variable $h\to 0$, under the correspondence $x=1/h$ the sectorial neighborhood of infinity in the direction $A^*$ becomes a sectorial neighborhood of the origin in the direction $A$, and we will talk about a resurgent function $g(h)$ (for $h\to 0$) in the direction $A$.

\subsubsection{Examples of resurgent functions.}
 
\paragraph{Example 1.} $h^\alpha$ 

The major corresponding to $h^{\nu}$ for $\nu\ne 1,2,3,..$ is 
$\frac{-1}{2i\sin(\pi\nu)} \frac{(-\xi)^{\nu-1}}{\Gamma(\nu)}$, and $\xi^{\nu-1}\frac{\Log \xi}{2\pi i \Gamma(\nu)}$ for $\nu=1,2,...$.

\paragraph{Example 2.} $\log h$.

\paragraph{Example 3.}  $\Phi(h):=e^{1/h}$, $\Phi(h):=e^{-1/\sqrt{h}}$. (~\cite[R\'es II.3.4]{CNP})

\paragraph{Example 4.}  $e^{-1/h^2}$ is zero as a resurgent function on the arc  $(-\pi/4,\pi/4)$, but does not give a resurgent function on any larger arc  because there it is no longer bounded by any function $e^{c/|h|}$ for $c\in \R$.

\subsection{Decomposition theorem for a resurgent function} \label{DecompThm}

Making rigorous sense of a formula of the type
$$\phi(h) \ \sim \ \sum_k e^{-c_k/h}(a_{k,0}+a_{k,1}h+a_{k,2}h^2+...), \ \  \ \  h\to 0,$$
may be done by explicitly writing an estimate of an error that occurs if we truncate the $k$-th power series on the right at the $N_k$-th term. Resurgent analysis offers the following attractive alternative: if $\phi(h)$ is resurgent, then the numbers $c_k$ can be seen as locations of the first sheet singularities of the major ${\bf\Phi}(\xi)$, which is expressed in terms of a decomposition of ${\bf\Phi}(\xi)$ in a formal sum of microfunctions that we are going to discuss now. Microfunctions, or the singularities of ${\bf\Phi}(\xi)$ at $c_k$ are related to power series $a_{k,0}+a_{k,1}h+a_{k,2}h^2+...$ through the concept of Borel summation, see section \ref{BorelSum}.

\subsubsection{Microfunctions and resurgent symbols} \label{MicrofAndResSymb}

{\bf Definition.} (see ~\cite[Pr\'e II.1]{CNP}) A \underline{microfunction} at  $\omega\in\C$ in the direction ${\check A}\subset {\mathbb S}^1$ is the datum of a sectorial germ at $\omega$ in direction $\check A$ modulo holomorphic germs in $\omega$; the set of such microfunctions is denoted by  
$${\pmb{\cal C}}^\omega(A) =  {\cal O}^\omega({\check A}) / {\cal O}_\omega . $$
We will write ${\pmb{\cal C}}(A)$ to mean ${\pmb{\cal C}}^0(A)$.
A microfunction is said to be \underline{resurgent} if it has an endlessly continuable representative. The set of resurgent microfunctions at $\omega$ in direction $\check A$ is denoted in ~\cite[R\'es I.3.0, p.178]{CNP} by ${\pmb{\cal R}}^\omega (A)$.

\begin{figure}[h]\includegraphics{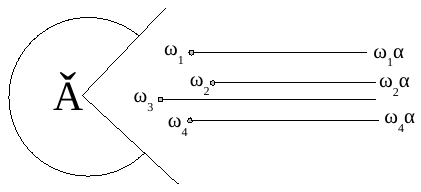} \caption{Singularities of a major and corresponding cuts} \label{ResDRWp8}
\end{figure}

{\bf Definition.} (~\cite[R\'es I.3.3, p.183]{CNP}) A \underline{resurgent symbol} in direction $\check A$ is a collection $\dot\phi = (\phi^\omega \in {\pmb{\cal R}}^\omega (A))_{\omega\in\C}$   such that $\phi^\omega$ is nonzero only for $\omega$ in a discrete subset $\Omega\in\C$ called the \underline{support} of $\dot\phi$, and for any $\alpha\in A$ the set $\C \backslash \Omega\alpha$ is a sectorial neighborhood of infinity in direction $\check A$, where  $\Omega\alpha=\bigcup_{\omega\in\Omega}\omega\alpha$ be the union of closed rays in the direction $\alpha$ emanating from the points of $\Omega$, Fig.\ref{ResDRWp8}. 

The set of such resurgent symbols is denoted $\dot{\pmb{\cal R}}(A)$, and resurgent symbols themselves can be written as  $\dot \phi = (\phi^\omega)_{\omega\in \Omega} \in \dot{\pmb {\cal R}}(A)$ or as $\dot \phi = \sum_{\omega\in\Omega} \phi^\omega \in \dot{\pmb {\cal R}}(A)$.

{\bf Definition.} A resurgent symbol is \underline{elementary} if its support $\Omega$ consists of one point. It is \underline{elementary simple} if that point is the origin. 

\subsubsection{Decomposition isomorphism.} \label{DecomposnIsomsm}

The correspondence between resurgent symbols in the direction $A$ and majors of resurgent functions in direction $A$ depends on a {\it resummation direction} $\alpha\in A$, which we will fix once and for all. The direction $\alpha$ can more concretely be thought of as $\arg h$ or as the the direction of the cuts in the $\xi$-plane 
that we are going to draw.

Suppose (~\cite[p.182-183]{CNP}) $\Omega$ is a discrete set (of singularities) in the complement of a sectorial neighborhood of infinity in direction $\check A$, and let a holomorphic function ${\bf \Phi}$ be defined on $\C\backslash \Omega\alpha$ and  endlessly continuable. Take  $\omega\in \Omega \alpha$.   Let $D_\omega$ be a small disc centered at $\omega$. Its diameter in the direction $\alpha$ cuts $D_\omega$ into the left and right open half-discs $D^{-}_\omega$ and $D^{+}_\omega$ (or top and bottom if $\alpha$ is the positive real direction). If $D_\omega$ is small enough, the function ${\bf\Phi}|D^{+}_\omega$, resp ${\bf\Phi}|D^{-}_\omega$, can be analytically continued to the whole split disc $D_\omega \backslash \omega \alpha$. Denote by $sing^{\omega}_{\alpha+} {\bf \Phi}$, resp. $sing^{\omega}_{\alpha-} {\bf \Phi}$, the microfunction  at $\omega$ of direction $\check \alpha$ defined by the class modulo ${\cal O}_\omega$ of this analytic continuation. 

\begin{Thm} (\cite[R\'es I.4]{CNP}) Let $\dot \phi = (\phi^\omega)_{\omega\in \Omega} \in \dot{\pmb{ \cal R}}(A)$ be a resurgent symbol with the singular support $\Omega$, and $\alpha \in A$. There is an endlessly continuable function ${\bf \Phi}\in {\cal O}(\C\backslash \Omega\alpha)$ such that
$$ sing^\omega_{\alpha+} {\bf \Phi} = \left\{ 
\begin{array}{ll} \phi^\omega & \text{if $\omega\in\Omega$,} \\ 0 & \text{if not.}
\end{array} \right. $$
This defines a bijective linear map
$$ {\bf s}_{\alpha_+} \ : \ \dot{\pmb{\cal R}}(A) \ \to \ {\pmb{\cal R}(A)}/{\cal O}(\C) \ : \  {\bf s}_{\alpha+}{\dot \phi} = {\bf \Phi} $$ whose inverse $({\bf s}_{\alpha+})^{-1}$ is called the \underline{decomposition isomorphism}. \end{Thm}

An endlessly continuable function ${\bf s}_{\alpha-}{\dot \phi}$ can be defined analogously. 

The maps ${\bf s}_{\alpha+}$,${\bf s}_{\alpha-}$  respect convolution products (cf. ~\cite[p.185, R\'es I.4]{CNP}). 

The map defined in ~\cite[R\'es I]{CNP} as
\begin{equation} ({\bf s}_{\alpha+})^{-1}\circ {\bf s}_{\alpha-} \ : \ {\dot {\pmb{\cal R}}}(A) \ \to \ {\dot{\pmb{\cal R}}}(A) \label{HomPassage} \end{equation}
is not identity and gives rise to the concept of alien derivatives which is central in resurgent analysis but will not be discussed in this article.

Further, if a resurgent function $\varphi(h)=\varphi(h,t)$ and its major ${\bf\Phi}(\xi)={\bf\Phi}(\xi,t)$ depend, say, continuously in some appropriate sense, on an auxiliary parameter $t$, the decompositon into microfunctions $({\bf s}_{\alpha+})^{-1}{\bf \Phi}$ will depend on $t$ ``discontinuously" -- an effect referred to as {\it Stokes phenomenon} and discussed, e.g., in ~\cite{DP99}.

\subsubsection{Mittag-Leffler sum} \label{MLS}

The concept of a Mittag-Leffler sum formalizes the idea of an infinite sum of resurgent functions $\sum_j \varphi_j(h)$ where $\varphi_j(h)$ have smaller and smaller exponential type, e.g., $\varphi_j = O(e^{-c_j/h})$ for $c_j\to \infty$ as $j\to \infty$.


Rephrasing ~\cite[Pr\'e I.4.1]{CNP}, let $\Phi_j$, $j=1,2,...$ be endlessly continuable holomorphic functions (thought of as majors of $\varphi_j$), $\Phi_j\in{\cal O}(\Omega_j)$, where $\Omega_j$ are sectorial neighborhoods of infinity satisfying $\Omega_j\subset \Omega_{j+1}$ and $\bigcup_j \Omega_j = \C$. Suppose ${\cal S}_j$ together with the projection $p_j:{\cal S}_j\to \C$ and with the choice of the first sheet (which contains $\Omega_j$), is the endlessly continuable Riemann surface of $\Phi_j$. The following statement seems to be implicitly used in ~\cite{CNP}; we will use it, too, although we are not aware of a detailed treatment of this question in the literature:

\begin{Statement} \label{RiemSurfMLS} Under above assumptions, there exists an endlessly continuable Riemann surface ${\cal S}_{ML}$ together with the choice of the first sheet and locally biholomorphic maps $\pi: {\cal S}_{ML} \to \C$ and $\pi_j:{\cal S}_{ML}\to {\cal S}_j$, so that: \\
(i)  $\pi=p_j \circ \pi_j$ and $\pi_j$ maps the first sheet of ${\cal S}$ to the first sheet of ${\cal S}_j$, for all $j$; \\
(ii) for any point $\xi\in{\cal S}_{ML}$, there is an $N\in \N$ so that $\pi_j(\xi)\in \Omega_j$ (and hence is on  the first sheet of ${\cal S}_j$) for all $j\ge N$. 
\end{Statement}

Then, similarly to ~\cite{CNP}, take an exhaustive sequence of discs $D_n\subset \Omega_n$, $\bigcup_n D_n=\C$, and find entire functions $E_j$ so that $\sup_{\xi\in D_j} |\Phi_j(\xi) -E_j(\xi)| \le 1/j!$. Then, thanks to the condition (ii) of the Statement \ref{RiemSurfMLS}, the series $\sum_j \pi_j^* (\Phi_j-p_j^* E_j)$ is factorially convergent on compact subsets of ${\cal S}_{ML}$.

By $\MLS_j \Phi_j$ we will mean any such sum $\sum_j  \pi_j^* (\Phi_j-p_j^* E_j)$ whose terms are {\it factorially convergent} on compact subsets of ${\cal S}_{ML}$; being specific about the rate of convergence will be important later on. We will call $\Phi$ a \underline{Mittag-Leffler sum} of $\Phi_1,\Phi_2,...$ and write
$$ \Phi \ = \ \MLS_j \Phi_j. $$ 


\subsection{Borel summation. Resurgent asymptotic expansions.} \label{BorelSum}

{\bf Definition.}  A \underline{resurgent hyperasymptotic expansion}  is a  formal sum 
$$\sum_k e^{-c_k/h}(a_{k,0}+a_{k,1}h+a_{k,2}h^2+...),$$ where: \\
i) $c_k$ form a discrete subset in $\C$ in the complement to some sectorial neighborhood of infinity in direction $\check A$;\\
ii) the power series of every summand satisfies the Gevrey condition, and \\
iii) each infinite sum $a_{k,0}+a_{k,1}h+a_{k,2}h^2+...$ defines, under formal Borel transform 
$$ {\cal B} \ : \ e^{-c_k/h }h^\ell \ \mapsto \ (\xi-c_k)^{\ell-1} \frac{\log (\xi-c_k)}{2\pi i \Gamma(\ell)} \ \ \text{if} \ \ell\in \N,$$ 
$$ {\cal B} \ : \ e^{-c_k/h } \ \mapsto \ \frac{1}{2\pi i(\xi-c_k)}, $$ 
an endlessly continuable microfunction centered at $c_k$. 

The authors of ~\cite{CNP} denote by  $\dot{{\cal R}} (A)$ (regular, as opposed to the bold-faced, ${\cal R}$) the algebra of resurgent hyperasymptotic expansions.

\begin{figure}[h]\includegraphics{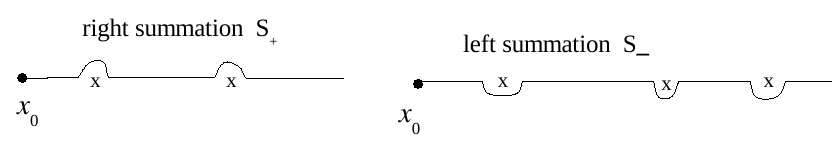} \caption{Integration contours in the definition of left and right summations.} \label{nuthesisp3}
\end{figure}

The right and left summations of resurgent asymptotic expansions are defined in ~\cite{DP99} or ~\cite{CNP} as follows. Given a Gevrey power series $\sum_{k=1}^{\infty} a_k h^k$, replace it by a function (the corresponding {\it``minor''}) ${\pmb f}(\xi)=\sum_{k=1}^{\infty} a_k \frac{\xi^{k-1}}{(k-1)!}$, assume that ${\pmb f}(\xi)$ has only a discrete set of singularities, and consider the Laplace integrals $\int_{[0,\alpha)}e^{-\xi/h}{\pmb f}(\xi)d\xi$ along a ray from $0$ to infinity in the direction $\alpha$ deformed to avoid the singularities from the right or from the left, as on the figure \ref{nuthesisp3}.   

After some technical discussion, this procedure defines a resurgent function of $h$ which ~\cite{CNP} denote 
${\rm S}_{\alpha \pm} \left(\sum_{k=0}^{\infty} a_k h^k \right)$.
Note that comparing the results of the left and right resummations is related the study of the map \eqref{HomPassage}.

\begin{figure}[h]
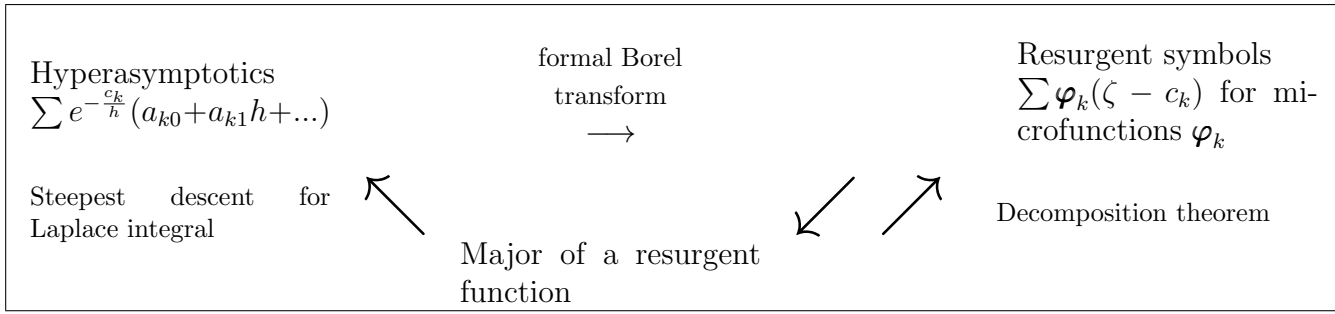

\fbox{\begin{tabular}{ccccc}
\parbox{4cm}{Hyperasymptotics
$\sum e^{-\frac{c_k}{h}}(a_{k0}+a_{k1}h+...)$} & & 
\parbox{3.5cm}{\begin{center} {\footnotesize formal Borel transform} \\  $\longrightarrow$ \end{center}}
&&
\ \ \ \  \parbox{4cm}{Resurgent symbols \\ 
$\sum {\pmb \varphi}_k(\zeta-c_k) $
for microfunctions ${\pmb\varphi}_k$} 
\\
 \parbox{4cm}{\footnotesize Steepest descent for Laplace integral} & {\huge $\nwarrow$}  && {\huge $\swarrow \ \nearrow$} & \parbox{4cm}{\footnotesize Decomposition theorem} 
\\
&& \parbox{4cm}{Major of a resurgent function} 
\end{tabular}} 
\caption{Logical relationship between concepts of resurgent analysis} \label{Diagram}
\end{figure}

The diagram on figure \ref{Diagram} helps  to visualize the logical relationship of the concepts that have been introduced. 

\section{Majors exponentially decreasing along a path.} \label{ExpDecreaseAlongPath}

As we have seen, the definition of the Laplace isomorphism ${\cal L}$ involves choosing a representative of a class $\mod \O(\C)$ of a major that is bounded along infinite branches of a contour. A result on the existence of such a  representative is recalled in section \ref{OneMajorExpDecr}; it might be helpful to look at the proofs in ~\cite[Pr\'e I.3]{CNP} before proceding to a generalization of this result to the case of an infinite sum of majors presented in section \ref{CaseOfConvergSeriesOfMajors}.

We write $D_R=\{ z\in \C : |z|\le R\}$.

\subsection{Case of a single major.} \label{OneMajorExpDecr}

\begin{Lemma} (\cite[Pr\'e I.3, Lemma 3.0]{CNP})  Let $\Gamma\subset \C$ be an embedded curve (i.e. a closed submanifold of dimension 1) transverse to circles $|\zeta|=R$ for all  $R\ge R_0$. Let ${\bf \Phi}$ be a holomorphic function in a neighborhood of $\Gamma$. Then for any  function $m:\R_+\to \R_+$ satisfying $\inf_{x\le N} m(x) >0$ for any $N>0$, there is an entire function $E$ such that $|({\bf \Phi} + {\bf E})(\zeta)|\le m(|\zeta|) $ for any $\zeta\in\Gamma$. \label{Lemma21} \end{Lemma}





The following two auxiliary results used by ~\cite{CNP} in the proof of lemma \ref{Lemma21} will be needed later on.

\begin{figure}[h]\includegraphics{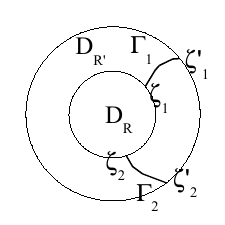} \caption{Notation in Lemma \ref{L321}.} \label{ResDRWp105}
\end{figure}

\begin{Lemma} \label{L321} Let (cf. Fig.\ref{ResDRWp105}) $\Gamma\subset \C$ be a finite disjoint union of curvilinear segments $\Gamma_i = [\zeta_i,\zeta'_i]$ with $i=1,...,r$, joining transversally the boundaries of the annulus $R\le |\zeta| \le R'$. 
Let ${\bf \Phi}$ be a function holomorphic in a neighborhood of $\Gamma$ and zero at $\zeta_1,...,\zeta_r$.
 Then for any $\varepsilon>0$ there is a polynomial function ${\bf E}$ such that\\
i) $|{\bf E}|\le \varepsilon$ on $D_R = \{ |\zeta|<R\}$; \\
ii) $|{\bf \Phi}+{\bf E}|\le \varepsilon$ on $\Gamma$; \\
iii) $({\bf \Phi}+{\bf E})(\zeta'_i)=0$ for $i=1,...,r$. \end{Lemma}

\begin{Lemma} Given $r$ points $\zeta_1,...,\zeta_r$ on a compact set $K\subset \C$, then for all $\varepsilon>0$ there is $\varepsilon''>0$ such that for any choice of the interpolation data $a_1,...,a_r$ satisfying $|a_i|\le\varepsilon''$, the Lagrange interpolation polynomial defined by $Q(\zeta_i)=a_i$ is estimated by $\varepsilon$ on $K$.  \end{Lemma}



\subsection{Case of a convergent series of majors.} \label{CaseOfConvergSeriesOfMajors}

In order to be able to work with a Laplace integral of an infinite sum of majors, it is useful to choose a representative of each summand of that series that would be bounded on the infinite branches of the integration path, and to do it in a way that such choices would be consistent with taking the infinite sum. The main issue is to show that the entire correction functions to each of the summands form themselves a series convergent on compact subsets of $\C$. More precisely, 

\begin{Prop} \label{Prop32} Suppose $\sum_{j=1}^{\infty} \Phi^{(j)}$ is an infinite series of majors defined on the same (unbounded) domain ${\cal S}\subset \C$ which converges uniformly and faster than some geometric series with a ratio $q$ (i.e., for every disc $D_\ell$ there is a constant $M_\ell$ so that $|\Phi^{(j)}|\le M_\ell q^j$ on $D_\ell\cap {\cal S}$) 
on compact subsets to a function $\Phi$,  $\Gamma$ a contour transversal to circles $\partial D_R$ for $R\ge R_0$  and $m:\R_+\to \R_+$ a function as in Lemma \ref{Lemma21}.  Then for any number $Q$, $q<Q<1$, one can choose  entire functions $F^{(j)}$ so that: \\
i) $\sum_{j=1}^{\infty}  F^{(j)}$ converges uniformly on compact sets  of $\C$ faster than geometric series of ratio $Q$ to an entire function $F$;\\
ii) for all $j$ the function $|\Phi^{(j)}(\xi) +  F^{(j)}(\xi)|\le (1-q) Q^j m(|\xi|)$  along $\Gamma$;\\
iii) the function $|\Phi(\xi)+ F(\xi)|\le \frac{1-q}{1-Q}m(|\xi|)$ along $\Gamma$; \\
iv) $|\Phi^{(j)} + F^{(j)}|$ can be, on every compact subset of ${\cal S}$, estimated by a geometric series with the ratio $Q$. \end{Prop}

\textsc{Proof.} Only the parts i) and ii) need to be proven in detail, the parts iii) and iv) will then follow as easy consequences.

Suppose, to simplify notation, that  $R_0=1$, that $m(|\zeta|)\le 2^{-|\zeta|}$ and $m(|\zeta|)=m_k$ (constant) for $k-1<|\zeta|\le k$. 

There are positive numbers $b_{k\ell}$ satisfying the following property: \\ 
{\it For any interpolation data on the the finite set  $\Gamma\cap \partial D_k$, namely, any function $d:\Gamma\cap \partial D_k \to \C$ with $\max_{p\in \Gamma\cap \partial D_k}|d(p)|\le 1$, the supremum of the corresponding degree $|\Gamma\cap \partial D_k|-1$  interpolation polynomial is $<b_{k\ell}$ on $D_\ell$.} \\
 Clearly we can choose $b_{k\ell}\le b_{kk}$ for $\ell\le k$ and $b_{k\ell}\ge 1$ for all $k,\ell$. 

Fix numbers  $\{Q_s\}_{s\in\N}$ such that $q<Q_1<Q_2<...<Q_s<...<Q<1$.

We will construct by induction entire functions $F^{(j)}_1,..., F^{(j)}_{s}$ such that for each fixed $s$, the series $\sum_{j=1}^\infty [\Phi^{(j)}+F_1^{(j)}+...+F^{(j)}_{s}]$ converges uniformly on compact sets faster than a geometric series of ratio $Q_{s}$.

\underline{For $s=1$},
find a number $N_1$ so that for $\forall j> N_1$ one has $|\Phi^{(j)}(\xi)| < \frac{(1-Q_1)Q_1^j}{2b_{11}} m_1$ on $\Gamma\cap D_1$.
Using lemma \ref{Lemma21}, choose entire functions $E^{(j)}_1$ so that $|\Phi^{(j)}+E^{(j)}_1|<\frac{m_1(1-Q_1)Q_1^j}{b_{11} 2}$ on $D_1\cap \Gamma$ for $j\le N_1$ and $E^{(1)}_j=0$ for $j> N_1$.
Now choose $G^{(j)}_1$ as interpolation polynomials so that $\Phi^{(j)}+E^{(j)}_1+G^{(j)}_1=0$ on $\Gamma \cap \partial D_1$. Then 
\begin{equation} \sup_{D_\ell} |G^{(j)}_1| < \frac{b_{1\ell} m_1(1-Q_1)Q_1^j}{2b_{11}} \ \ \ \text{for all} \ j\in \N . \label{eq818} \end{equation} 

Put $F^{(j)}_1:=E^{(j)}_1+G^{(j)}_1$. Then on $D_1\cap \Gamma$
\begin{equation}  | \Phi^{(j)} + F^{(j)}_1 | \ \le \ |\Phi^{(j)}+E^{(j)}_1| + |G^{(j)}_1| \le \frac{m_1 (1-Q_1)Q_1^j}{2b_{11}} + \frac{m_1(1-Q_1)Q_1^j}{2}\ \le\ m_1(1-Q_1)Q_1^j \label{eq241} \end{equation}

On the other hand, on $D_\ell$ for $j\ge N_1 +1$,
$$  |\Phi^{(j)}+F_1^{(j)}| =  |\Phi^{(j)}+G^{(j)}| \le |\Phi^{(j)}|+|G^{(j)}|, $$
where the first summand decreases faster then some geometric series of ratio $Q_1$ by assumptions of the proposition, and the second does so by \eqref{eq818}.

\underline{For $s\ge 2$,} suppose that we have constructed $F^{(j)}_1,...,F^{(j)}_{s-1}$ for all $j$, let us construct $F^{(j)}_s$. Choose $N_s$ so that on $\Gamma \cap (D_s\backslash D_{s-1})$ we have $\forall j>N_s$ the inequality  $|\Phi^{(j)}+F_1^{(j)}+...+F_{s-1}^{(j)}| < \frac{m_s(1-Q_s)Q_s^j}{2b_{ss}}$. Then, by lemma \ref{L321}, there are entire functions $E^{(j)}_s$ such that\\ 
a) 
\begin{equation}  |E^{(j)}_s|<\frac{m_s(1-Q_s)Q_s^j}{2b_{ss}}\ \ \text{on} \ D_{s-1} \ \ \text{for} \ j\le N_s \label{eq900} \end{equation}
and \\
b) $|E^{(j)}_s+F^{(j)}_1+...+F_{s-1}^{(j)} + \Phi^{(j)}|<\frac{m_s(1-Q_s)Q_s^j}{2b_{ss}}$ on $(D_s\backslash D_{s-1}) \cap \Gamma$. \\
Choose $G_s^{(j)}$ as the interpolation polynomial such that $F^{(j)}_{s-1}+...+ F^{(j)}_1+ \Phi^{(j)}+E^{(j)}_s+G^{(j)}_s = 0$ on $\Gamma \cap \partial D_s$ and put $F^{(j)}_s:= E^{(j)}_s+G^{(j)}_s$. Then \begin{equation} \sup_{D_\ell} |G_s^{(j)}| < \frac{b_{s\ell} m_s (1-Q_s)Q_s^j}{2b_{ss}} \ \ \text{ for all } \ j \label{eq819} \end{equation}  

Combining \eqref{eq900} and \eqref{eq819}, we see  that
\begin{equation}  \sup_{D_{s-1}} |F^{(j)}_s|  \le m_s (1-Q_s) Q_s^j \label{eq902} \end{equation} 

Analogously to the \eqref{eq241},
\begin{equation} |\Phi^{(j)}+F_1^{(j)}+...+F_s^{(j)}|\le m_s(1-Q_s)Q_s^j  \label{eq242} \end{equation}
on $\Gamma\cap (D_s\backslash D_{s-1}).$

The series $\sum_j (\Phi^{(j)}+F_1^{(j)}+...+F_s^{(j)})$ still converges uniformly and faster than some geometric series of ratio $Q_s$ on compact sets; indeed, on $D_\ell$ such an estimate follows from the analogous property for $(s-1)$, the fact that for $j$ large $F^{(j)}_s=G^{(j)}_s$ , and the formula \eqref{eq902}.

This finishes the inductive construction. 

Now put $ F^{(j)} := \sum_{k=1}^\infty F^{(j)}_k $. 
This sum is uniformly convergent on compact sets because of the estimate $|F^{(j)}_k|\le m_k\le \frac{1}{2^k}$ that holds on $D_\ell$ for $k\ge \ell+1$ by \eqref{eq902}. Hence $F^{(j)}$ is a well-defined entire function.


Let us now find a geometric series with ratio $Q$ that bounds $F^{(j)}$ on $D_\ell$, at least for $j>\max\{ N_1,...,N_\ell\}$. We get
$$ \sup_{D_\ell} |F^{(j)}| \ \le \  \sup_{D_\ell} \sum_{k=1}^\infty |F^{(j)}_k| \ \le \ \sum_{k=1}^{\ell} \sup_{D_\ell} |F^{(j)}_k| + \sum_{k=\ell+1}^{\infty} \sup_{D_\ell} |F^{(j)}_k| \ \stackrel{(!)}{\le} \ $$
$$ \ \le \ \sum_{k=1}^{\ell} b_{k\ell}m_k (1-q) Q^j + \sum_{k=\ell+1}^{\infty}  m_k (1-q) Q^j \ \le \ 
\left[ \sum_{k=1}^{\ell} b_{k\ell}m_k (1-q)  + \sum_{k=\ell+1}^{\infty}  m_k (1-q) \right] Q^j, $$
where the inequality (!) holds because in the first sum in the given range we can use the equality $F^{(j)}_k=G^{(j)}_k$ and \eqref{eq819}, and for the second sum we use can \eqref{eq902}. This proves that $F^{(j)}$ satisfies (i). 

The statement (ii) follows from \eqref{eq241} and \eqref{eq242}. $\Box$


\subsection{Interchangeability of infinite sum and ${\cal L}$.}

\begin{Prop} \label{Prop35}  Suppose $\Phi^{(j)}$ is a sequence of majors defined on a common sectorial neighborhood of infinity and suppose that on each compact subset of this neighborhood the sum $\sum_j |\Phi^{(j)}|$ uniformly converges  faster than a geometric series with a ratio $q<1$. Then 
$${\cal L} \left\{ \left( \sum_j \Phi^{(j)}\right) \mod {\cal O}(\C) \right\} \ = \ \sum_j {\cal L} (\Phi^{(j)} \mod {\cal O}(\C)).$$ 
\end{Prop}

\textsc{Proof.} Taking a number $Q$, $q<Q<1$, a path $\gamma$ as in \eqref{defLfla}, the representatives $\Phi^{(j)}$ of the corresponding integrality classes and the number $Q<1$   provided by proposition \ref{Prop32} for $m(t)=2^{-t}$, we can write  ${\cal L}(\Phi^{(j)} \mod {\cal O}(\C))$ as $\int_\gamma e^{-\xi/h} \Phi^{(j)} d\xi$, and similarly for ${\cal L}\left\{ \left( \sum_j \Phi^{(j)}\right) \mod {\cal O}(\C) \right\}$. Then the question reduces to showing that
$$ \int_\gamma \sum_j e^{-\xi/h} \Phi^{(j)}(\xi) d\xi \ = \ \sum_j \int_\gamma  e^{-\xi/h} \Phi^{(j)}(\xi) d\xi.$$
By Fubini's theorem we need to check:\\

i) For any $\xi\in \gamma$ the sum $\sum_j e^{-\xi/h} |\Phi^{(j)}(\xi)|$ converges -- because  $|\Phi^{(j)}(\xi)|<m(|\xi|) Q^j$; \\
ii) the integral of such a sum clearly converges for small $h$; \\
iii) For any $j$ the integral $\int_\gamma  e^{-\xi/h} |\Phi^{(j)}(\xi)| d\xi$ converges -- in fact it is less than $\frac{1}{Q^j} \int_\gamma |e^{-\xi/h}| d|\xi|$; \\
iv) the sum of these integrals is then clearly convergent, too. $\Box$

\begin{Cor} \label{CorLMLS} Given a sequence of majors $\Psi_j$ as in section \ref{MLS}, then
$$ {\cal L} \left( \MLS_j \Psi_j \right) =  \sum_j {\cal L} \Psi_j $$ $\Box$ \end{Cor}


\section{Convolution products and majors bounded in a neighborhood of infinity.} \label{ConvolutionSectn}

Let $A$ be a small arc in the circle of directions which for simplicity of language will be assumed symmetric with respect to the real axis. Let us study the convolution product of majors that will be assumed holomorphic on some sectorial neighborhood of infinity in the direction ${\check A}$. This convolution of majors is known to correspond to multiplication of resurgent functions.

Recall ~\cite{CNP} that the convolution of two integrality classes of majors $[\Phi]$ and $[\Psi]$ along a path $\Gamma$ adopted to a sectorial neighborhood of infinity in direction $\check A$ is defined by choosing two representatives: $\Phi$ that is exponentially decreasing along the infinite branches of $\Gamma$, and $\Psi$ that is $\le const$ on the neighborhood of infinity bounded by $\Gamma$ (see lemma \ref{PreI342} below), and considering the integral
 $$ (\Phi*_\Gamma \Psi)(\xi) = \int_\Gamma  \Phi(\eta)\Psi(\xi-\eta) d\eta. $$
This integral defines a sectorial germ of an analytic function at infinity and can be analytically continued to some Riemann surface by deforming the integration contour.
 
It is discussed in ~\cite{CNP} that the result of the convolution is independent of choices modulo ${\cal O}(\C)$.

\begin{figure}[h]\includegraphics{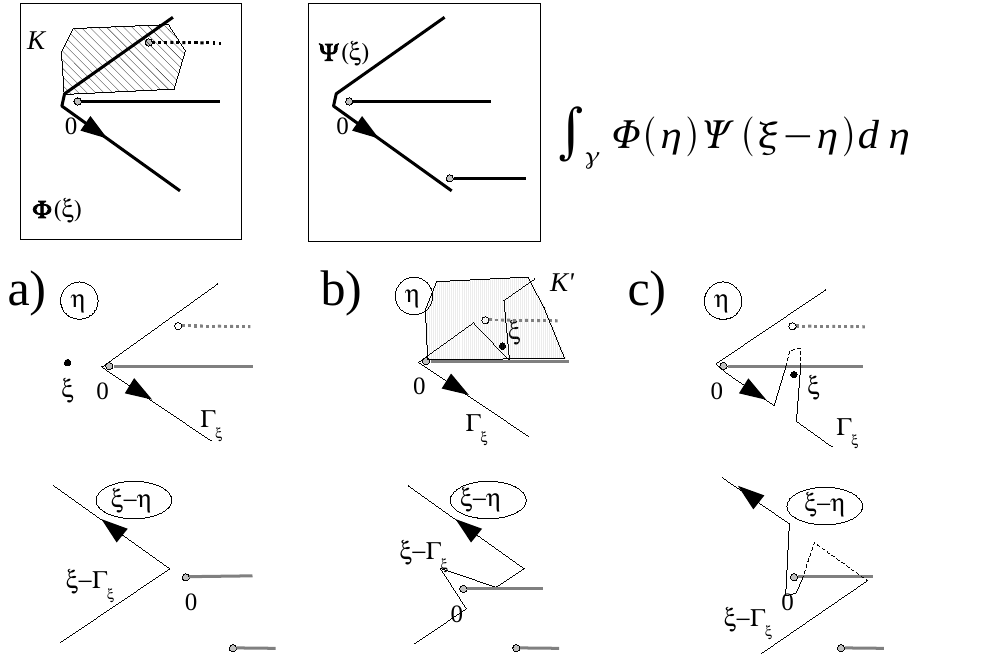} \caption{Deformation of the convolution contour and analytic continuation of the convolution product.} \label{ExNote3p3}
\end{figure}

Let us look a little closer at the deformation of the contour. Suppose for simplicity that $\Gamma$ consists of two rays starting at $\eta=0$. Then for $\xi$ to the left of $\Gamma$ the above formula defines an analytic function ``on the nose"; this is the case a) on the figure \ref{ExNote3p3}.

In order to analytically continue the convolution to other values of $\xi$, to the right of the contour $\Gamma$, we need to continuously deform the convolution contour (cases b), c) of the fugure \ref{ExNote3p3}) to obtain a contour $\Gamma_\xi$ so that $\Gamma_\xi$ avoids singularities of $\Phi$ and $\xi-\Gamma_\xi$ avoids singularities of $\Psi$. The singularities of the convolution appear for those $\xi$ for which this deformation becomes impossible.

It is shown in ~\cite{CNP} that if $\Phi$ and $\Psi$ are defined on endless Riemann surfaces ${\cal S}_\Phi$ and ${\cal S}_\Psi$ respectively, then $\Phi*\Psi$ can be analytically continued to an endless Riemann surface denoted ${\cal S}_\Phi * {\cal S}_\Psi$.

Note that if $\xi$ stays within a compact subset $K$ ($K$ is a subset of ${\cal S}_\Phi *{\cal S}_\Psi$, but we thought it helpful to superimpose it on the Riemann surface ${\cal S}_\Phi$ on the picture), then the deformation of the contour $\Gamma_\xi$ will be confined to a compact subset $K'$ of ${\cal S}_\Phi$.  

As has been mentioned, the definition of the convolution product involves a choice in an integrality class of a major, $\Psi \mod \O(\C)$, of a representative that is bounded on sectorial neighborhood of infinity in the direction $\check A$. Let us begin by recalling how this choice is made for a single major, and then prove that this choice can be made compatible with an infinte sum of majors (section \ref{CaseOfConvergSeriesOfMajors42}). Finally, we will show that under appropriate assumptions convolution is interchangeable with an infinite sum of majors (section \ref{InfiniteSumConvoRev}).

\subsection{Case of a single major} \label{CaseSingleM}

\begin{Lemma} \label{PreI342} (cf. ~\cite[Pr\'e I.3.4.2]{CNP})  Every integrality class $[\Phi]$ in the direction $\check A$ for a small  arc $A$ has a representative $\Psi$  bounded in a (possibly smaller) sectorial neighborhood in the direction $\check A$. \end{Lemma}

\begin{figure}[h]\includegraphics{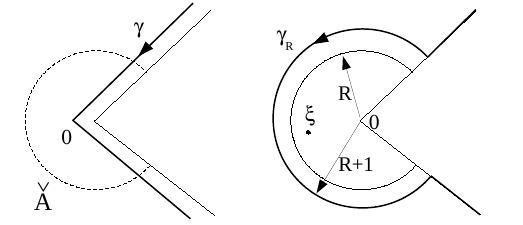} \caption{Integration contours in the proof of Lemma \ref{PreI342}} \label{ExNote3p2}
\end{figure}

\textsc{Proof.} Let $\Omega$ be a sectorial neighborhood of infinity in the direction ${\check A}$ and suppose $\Phi$ is defined in $\Omega$. Let $\gamma$ be an infinite path contained in $\Omega$ and adapted to ${\check A}$ consisting of two rays coming together, for simplicity of notation, at the origin, fig.\ref{ExNote3p2}. Without loss of generality assume $\Phi$ to be exponentially decreasing along the branches of $\gamma$. 

Then
$$\Phi(\xi)  \ = \ \Psi(\xi)  \ + \  E(\xi),$$
where
$$ \Psi(\xi) \  = \frac{1}{2\pi i} \int_\gamma \frac{\Phi(\eta)}{\xi-\eta} d\eta $$
and where $E\in {\cal O}(\C)$ can be defined as follows: \\
Given $\xi\in D_R\subset \C$, construct the contour $\gamma_R$ consisting of an arc of radius $R+1$ and two infinite branches of $\gamma$ as shown on the fig.\ref{ExNote3p2} and put
$$ E(\xi) = \frac{1}{2\pi i}\int_{\gamma_R} \frac{\Phi(\eta)}{\xi-\eta}d\eta. $$

Finally, notice that $\Psi$ is bounded in any subset of $\Omega$ where $\sup_{\eta\in\gamma} \frac{1}{|\xi-\eta|}$ is bounded.  $\Box$

We have recalled this proof in order to make the following

{\bf Remark \ref{PreI342}.A.} If $\xi\in D_R$,  and if $\Phi$ is bounded by an exponentially decreasing function $m(|\xi|)$  on $\gamma$ and by $M_{R+1}$ on ${\rm Sec}(0,\check A)\cap D_{R+1}$, we have an estimate 
$$ |E(\xi)|  \le  \frac{1}{2\pi} \int_{arc} |\Phi(\eta) | |d\eta|   +  \int_{\gamma} m(|\xi|) |d\xi| \le (R+1) M_R  +  \int_{\gamma} m(|\xi|) |d\xi|. $$  

That means that if $\Phi^{(j)}$ are bounded by $q^j m(|\xi|)$ on the infinite branches of $\gamma$ and converge faster than some geometric series with ratio $q$ on $D_{R+1}$, then on $D_R$ the values of $E(\xi)$ are bounded by some geometric series  with  ratio $q$.

\subsection{Case of a convergent series of majors} \label{CaseOfConvergSeriesOfMajors42}

The following proposition will be used in the proof of Prop.\ref{Prop9} later.


\begin{Prop} \label{Prop331} Given a series of majors $\Phi^{(j)}$, all of them analytic on  ${\rm Sec}(0,\check A)$, converging uniformly and faster than a geometric series of ratio $q<1$ on compact subsets of their common Riemann surface ${\cal S}$. Let $\gamma$ be a contour contained in a sector ${\rm Sec}(p_0,\check A)$ and adapted to ${\check A}$, and let $p \in {\rm Sec}(p_0,\check A)$ be such that the distance from ${\rm Sec}(p,\check A)$  to $\gamma$ is $\varepsilon>0$. Then for any number $Q$, $q<Q<1$, we can choose entire functions $F^{(j)}$  so that: \\
i) In ${\rm Sec}(p,\check A)$ and on the compact subsets $K$ of ${\cal S}$, $|\Phi^{(j)}-F^{(j)}|<M_K Q^j$. \\
ii) $\sum_j  F^{(j)}$ converges uniformly on compact subsets.
\end{Prop}

\textsc{Proof.}  
Indeed, begin by choosing $\Phi_1^{(j)}$ such that $\Phi_1^{(j)}(\xi)\le Q^j m(|\xi|) $ along $\gamma$ with $m(|\xi|)=e^{-|\xi|}$ and such that  $\Phi^{(j)}-\Phi_1^{(j)}$ are holomorphic and converge faster than geometric series of ratio $Q$ on compact subsets.

Choose $E^{(j)}$ similarly to the proof of lemma \ref{PreI342}, 
$$E^{(j)}:= -\frac{1}{2\pi i} \int_\gamma \frac{\Phi_1^{(j)}(\eta)d\eta}{\xi-\eta} \ + \ \Phi_1^{(j)}(\xi); $$
then by remark \ref{PreI342}.A, $\sum_j E^{(j)}$ converges as geometric series of ratio $Q$ on compact subsets of $\C$ and hence defines an entire function. Take $F^{(j)} = (\Phi^{(j)}-\Phi^{(j)}_1)+E^{(j)}$.


Then we need to show that $|\Phi_1^{(j)}-E^{(j)}|<{\tilde M}_K Q^j$ on a compact subset $K$. ( Since, as noted before, a similar inequality holds for $|\Phi^{(j)}-\Phi_1^{(j)}|$, the differences $|\Phi^{(j)}-E^{(j)}|$ will also be estimated by a geometric series of ratio $Q$.) I.e., we need to show that $\int_\gamma \frac{\Phi_1^{(j)}(\eta)}{\xi-\eta}d\eta$ with $\gamma$ chosen as in the statement,  converges uniformly of compact subsets.

 Without loss of generality suppose $p=0$ and let $U_-={\rm Sec}(0,\check A)$ on the first sheet of ${\cal S}$, $U_+={\cal S}\backslash U_-$. For a compact subset $K\subset {\cal S}$ let $K_{\pm}=K\cap U_{\pm}$.


Fix $K\subset\subset {\cal S}$ and let us check  for $\xi\in K$ the inequality $\left| \int_\gamma \frac{\Phi_1^{(n)}(\eta)}{\xi-\eta} d\eta \right| < C Q^n$.

Suppose first that $\xi\in U_{-}$. We know that $|\Phi_1^{(n)}(\xi)|<C_1 Q^n e^{-|\xi|}$ on $\gamma$. When $\xi\in U_-$ and $\eta\in \gamma$, then $\xi-\eta\in U_-$ and so $\frac{1}{\xi-\eta}< C_2$. Then $\left| \int \frac{\Phi_1^{(n)}(\eta)}{\xi-\eta} d\eta \right| <C_3 Q^n$ for yet another constant $C_3$.   

Suppose now $\xi\in K_+$. In this case $\gamma$ gets deformed to a path $\gamma_\xi$ so that both $\gamma_\xi$ and $\xi-\eta$ are fully contained in $K'\cup U_-$ for some compact set $K'$ and the length of $\gamma_\xi \cap K'$ is $\le L_K$. Use the following estimates: for $\eta\in K'$ we have $|\frac{1}{\xi-\eta}|<C_4$ (note that $\xi$ is never equal to $\eta$ because $\xi\not\in \gamma_\xi$), $|\Phi_1^{(n)}(\eta)|<C_5 Q^n$ (because $\Phi_1^{(n)}$ converges faster than some geometric series on compact subsets, cf. Proposition \ref{Prop32},iv) \ ), so $\left| \int_{\gamma_\xi \cap K'} \frac{\Phi_1^{(n)}(\eta)}{\xi-\eta} d\eta \right|< L_K C_4 C_5 Q^n$.

For the part of the integral along $\gamma_\xi \backslash K'$ proceed as in the case of $\xi\in U_-$. $\Box$


\subsection{Interchanging infinite sum and a convolution.} \label{InfiniteSumConvoRev}


\begin{Prop} \label{Prop43rev} Suppose $\Psi^{(j)}$ is a sequence of majors defined on a common Riemann surface containing ${\rm Sec}(p_0,\check A)$ and converging faster than a geometric series on compact subsets. Let $p\in {\rm Int} ({\rm Sec}(p_0,\check A))$ and $\Gamma$ be the contour consisting of two rays on the boundary of ${\rm Sec}(p,\check A)$.  Then 
$$ \Phi *_{\Gamma} \left( \sum_j \Psi^{(j)} \right) \ = \ \sum_j \Phi *_\Gamma \Psi^{(j)} \ \ \ \ \mod {\cal O}(\C). $$
Moreover, with the appropriate choice of the representatives, the series in the right hand side converges uniformly on compact sets faster than a geometric series of some ratio $Q<1$.
\end{Prop}

\textsc{Proof.} Without loss of generality $p=0$; denote $U_- ={\rm Sec}(0,\check A)$.
We are studying the integrals $\int_\Gamma \Phi(\xi-\eta)\Psi^{(n)}(\eta) d\eta$. Choose representatives  $\Psi_n$ of $\Psi^{(n)}$ provided by Proposition \ref{Prop32} such that for some $Q<1$ we have $|\Psi^{(n)}(\xi)|\le M e^{-|\xi|} Q^n$ along $\Gamma$. 
Choose a representative $\Phi$ so that $|\Phi(\xi)|<C_1 $ on $U_-\cup \Gamma$. 

Let us show that with this choice of representatives the equality holds exactly, not just modulo entire functions. By Fubini's theorem, we need to show that $\int$, $\int\sum$, $\sum$, $\sum\int$ are absolutely convergent. Let us show that sum of analytic continuations of $\int_\Gamma \Phi(\xi-\eta)  \Psi_n (\eta) d\eta $ is convergent uniformly on compact sets of the Riemann surface ${\cal T} = {\cal S}_{\Phi} * {\cal S}_\Psi$. 

 Let us show that $\sum_{n=0}^\infty  \int_\Gamma \Psi_n(\eta)\Phi(\xi-\eta)d\eta$ converges uniformly on compact subsets. Fix $K\subset {\cal T}$ (${\cal T}={\cal S}_\Psi *{\cal S}_{\Phi^{(j)}}$ and let us check  for $\xi\in K$ the inequality $\left| \int_\Gamma \Psi(\eta)\Phi_n(\xi-\eta)d\eta \right| < C\alpha^n$.  For our compact set $K$, let $K_-=K\cap U_-$.


Suppose first that $\xi\in U_{-}$.  When $\xi\in U_-$ and $\eta\in \Gamma$, then $\xi-\eta\in U_-$; for $\eta\in\Gamma$ we also have  $|\Psi_n(\eta)| < C_2 Q^n e^{-|\eta|}$. In this case we can take $\Gamma_\xi=\Gamma$ and then $|\int_{\Gamma_\xi} \Phi(\xi-\eta) \Psi_n(\eta) d\eta | <C_3 Q^n$ for yet another constant $C_3$.

Suppose now $\xi\in K_+ := K\backslash K_-$.  In this case $\Gamma$ gets deformed to a path $\Gamma_\xi$ that both $\Gamma_\xi$ and $\xi-\eta$ are fully contained in $K'\cup U_-$  for some compact set $K'$ and the length of $\Gamma_\xi \cap K'$ is $\le L_K$. We will use the following estimates: for $\eta\in K'$ we have $|\Psi_n(\eta)|<C_4Q^n$, $|\Phi(\xi-\eta)|<C_5$, so $\left| \int_{\Gamma_\xi \cap K'} \Psi(\eta) \Phi_n(\xi-\eta) d\eta \right|< L_K C_4 C_5 Q^n$.

For the part of the integral along $\Gamma_\xi \backslash K'$ proceed as in the case of $\xi\in U_-$; together this will dispose of the case $\xi\in K_+$.
$\Box$

We thank the anonymous referee for a suggestion that led to streamlining of this proof.



\section{Interchanging infinite sums and the reconstruction isomorphism.}

In \ref{DecompThm} we have reminded the correspondence between resurgent symbols and endlessly continuable majors. This correspondence respects infinite sums in the following sense.

\begin{Prop} \label{Prop341} (Inspired by ~\cite[R\'es 3.2.5]{CNP}) Given an infinite series of resurgent microfunctions at $0$ whose representatives $\phi_j(\zeta)$ converge uniformly and faster than geometric series with ratio $q$ on compact subsets of their common Riemann surface. Then  for any $Q$, $q<Q<1$, and any direction $\alpha$ one can choose majors $\Psi_j \in {\bf s}_{\alpha+} [\phi_j] + {\cal O}(\C)$ defined on a common endlessly continuable Riemann surface ${\cal S}$, so that the series $\sum_j\Psi_j$ converges faster than a geometric series with ratio $Q$ on compact subsets of ${\cal S}$, and   ${\bf s}_{\alpha+} \sum_j \phi_j = \sum_j \Psi_j \ {\rm mod} \ {\cal O}(\C)$).
\end{Prop}

{\bf Proof.}  Choose $\Gamma$ as in the construction of ${\bf s}_{\alpha+}$ (~\cite[p.186]{CNP}) and choose, using Prop. \ref{Prop32},  representatives $\Phi_j$ of microfunctions $\phi_j$  satifsfying $\Phi_j(\zeta)< Q^j e^{-|\zeta|}$ ($Q<1$) along the infinite branches of $\Gamma$ and bounded by some geometric series of ratio $Q$ on every compact set of their common Riemann surface.

Then  also the integrals   $\Psi_j \ = \ \frac{1}{2\pi i} \int_\Gamma \frac{ \Phi_j(\eta)}{\xi-\eta} d\eta$ are bounded by some geometric series of ratio $Q$ on every compact set of ${\cal S}$.
This is shown by an obvious modification of the argument from sections \ref{CaseOfConvergSeriesOfMajors42}-\ref{InfiniteSumConvoRev}. 
$\Box$


\section{Small resurgent functions.} \label{SmallResFunSec}

{\bf Definition.} (~\cite[Pr\'e II.4, p.157]{CNP}) A microfunction $\varphi \in {\pmb{\cal C}}(A)$ is said to be a \underline{small microfunction } if it has a representative ${\bf\Phi}$ such that ${\bf\Phi} = o (\frac{1}{|\zeta|})$ uniformly in any sectorial neighborhood of direction $\check A$ for  $A' \subset\subset A$. \\ E.g.,  $h^\alpha$ for $\alpha>0$ satisfies that property.

The following definition has been somewhat modified compared to (\cite[R\'es II.3.2, p.219]{CNP}).

\begin{figure}[h]
\includegraphics{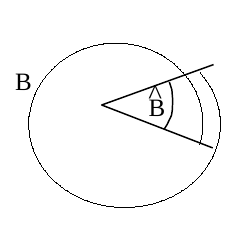} \caption{Arcs in the definition of a small resurgent microfunction.} \label{ExNote3p6}
\end{figure}

{\bf Definition.} For a given arc of direction $A=(\theta-\Delta\theta,\theta+\Delta\theta)$,  a \underline{small resurgent function} in the direction $A$ is such a resurgent function that all first-sheet singularities of one of (hence any of) its majors $\omega_\alpha$ satisfy $\re e^{-i\theta}\omega_\alpha > 0$, except maybe for one $\omega_0=0$,  and if $\omega_0=0$ then the corresponding microfunction is small in the direction of a large (i.e. $>2\pi$) arc $B$ with $\hat B \supset A$, see figure \ref{ExNote3p6}.

In the notation of this definition, we have:

\begin{Lemma} \label{MajorSmallResurgFctn} A small resurgent function in the direction $A$ can be represented by a major that is $o(1/|\xi|)$ in the direction $B'\subset\subset B$ around the origin.
\end{Lemma}

\textsc{Proof.} It is enough to prove the lemma for a small resurgent function whose decomposition consists of only one small resurgent microfunction at the origin, because the formal sum of all other microfunctions corresponds, via the decomposition isomorphism of section \ref{DecomposnIsomsm}, to a major that is holomorphic at the origin and is therefore automatically $o(1/|\zeta|)$ . This means that starting with a microfunction $[\phi]$ represented by an endlessly continuable sectorial germ $\phi\in {\cal O}^0(B)$, $\phi(\xi)=o(1/|\xi|)$ at the origin, we must construct a major $\Phi(\xi)$ of ${\bf s}_{\alpha+} [\phi]$ satisfying $\Phi(\xi)=o(1/|\xi|)$ for $\xi\to 0$ in the direction $B'$. Here $\alpha\in \hat B'$ will for definiteness be chosen as the positive real direction.

The major $\Phi$ will be given as an integral 
\begin{equation} \Phi(\xi)=\frac{1}{2\pi i}\int_{\Gamma_\xi} \frac{\phi(\eta)}{\xi-\eta}d\eta, \label{MajorSmallMinor}  \end{equation}
where the contour $\Gamma_\xi$ will be described presently. Let $\Delta\beta>0$ and $B'=\{ \theta: \ -\beta < \theta < 2\pi + \beta \} $, $B''=\{ \theta: \ -\beta-\Delta \beta < \theta < 2\pi + \beta +\Delta \beta \} $, $B=\{ \theta: \ -\beta-2\Delta \beta < \theta < 2\pi + \beta +2\Delta \beta \} $. Assume without loss of generality the length of $B$ to be $<\frac{5}{2}\pi$; otherwise cover $B$ by smaller arcs. Let $0<\beta'<\beta$ such that $\sin \beta'=\frac{1}{3}\sin (\beta+\Delta\beta)$. Choose a number $T>0$ so that $\phi(\xi)$ has no singularities for $\arg \xi\in B$ and $0<|\xi|<3T$;   and let $\xi$ belong to the sectorial neighborhood ${\cal U}$ of $0$ of radius $T$ in the direction $B'$.  (Remark that  $\Phi(\xi)$ can be analytically continued beyong ${\cal U}$ by deforming the contour, but the specific way of this deformation is not important at the moment.)

For $\xi\in {\cal U}$, the contour $\Gamma_\xi$ on the Riemann surface of $\phi$ will consist of a contour $\Gamma_\xi^0$ starting at $a_{2T}$ and ending at $b_{2T}$, 
 and the union  $\Gamma^{\infty}$ of two infinite branches independent of $\xi$, one from infinity to $a_{2T}$ and the other from $b_{2T}$ to infinity. Points $a_{2T}$ and $b_{2T}$ project to the point $2T$ on the complex plane, and to come from $a_{2T}$ to $b_{2T}$ one should go once counterclockwise around the origin.

\begin{figure}[h] \includegraphics{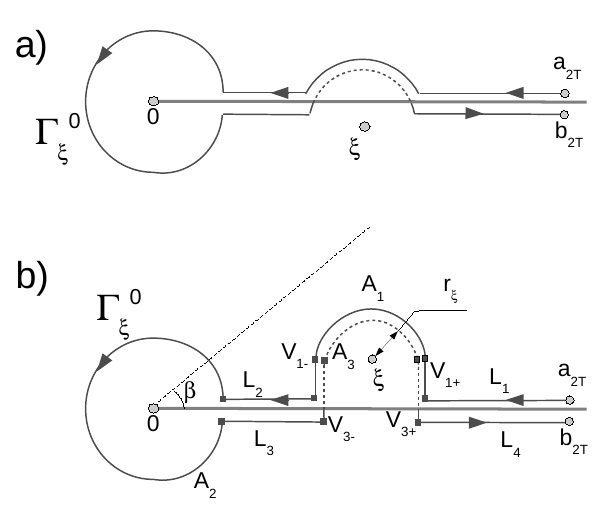} \ \ \ \parbox{12cm}{\small The inbound and outbound branches of the contour $\Gamma_\xi^0$ should be understood as lying on top of each other on different sheets of the Riemann surface of $\phi$. } \caption{Integration contour $\Gamma^0_\xi$ in Lemma \ref{MajorSmallResurgFctn}}
\label{Paper1p4}
\end{figure}

If $\beta \le \arg \xi \le 2\pi-\beta$, then the (projection of the) path $\Gamma_\xi^0$ will come from $\xi=2T$ along the positive real axis until point $|\xi|/3$, 
 go once counterclockwise around the circle $\{ \zeta \ : \ |\zeta|=|\xi|/3 \}$, and return along $[|\xi|/3, 2T]$. If $2\pi -\beta'  < \arg \xi \le 2\pi $,  then $\Gamma_\xi$ starts horizontally from $\xi=2T$, goes counterclockwise around the arc of the circle $S_\xi=\{ \zeta: | \zeta-\xi|=|\xi|\sin \beta' \}$,  back to the real line, one full circle counterclockwise around $\{ \zeta  :  |\zeta|=|\xi|/3 \}$, and retraces its path back to $2T$, fig.\ref{Paper1p4},a). 

If $2\pi \le \arg \xi \le 2\pi + \beta $, then instead of the arc along the circle $S_{\xi}$, the path contains two line segment $[\Re \xi \pm r_\xi; \xi \pm r_\xi]$ and the upper arc of the circle $|\zeta-\xi|=r_\xi$, where $r_\xi= \frac{1}{3}|\xi|\sin (2\pi+\beta+\Delta\beta-\arg \xi)$. 
Denote the horizontal parts of the path $L_1$, $L_2$, $L_3$, $L_4$,the vertical parts $V_{1\pm}$, $V_{3\pm}$ and the arcs $A_1,A_2,A_3$,  fig.\ref{Paper1p4},b).

For $-\beta<\arg \xi\le \beta$, $\Gamma^0_\xi$ can be defined similarly.

Two branches of $\Gamma^\infty$ on the Riemann surface of $\phi$ should be chosen to go in the postive real direction, slightly dispaced from the real axis so that their projections to the complex plane do not intersect except at $2T$, avoiding singularities of $\phi$ in some way, and satisfying conditions of Lemma \ref{Lemma21}. Then, by adding to it an entire function, we can assume $\phi(\xi)$ to be exponentially decreasing along $\Gamma^\infty$ and hence $\int_{\Gamma^\infty} \frac{\phi(\eta)}{\xi-\eta} d\eta$ is uniformly bounded for $\xi\in {\cal U}$.



We will prove the estimate $\Phi(\xi)=o(1/|\xi|)$ for $2\pi  \le \arg \xi \le 2\pi + \beta$, other cases are similar but simpler.

Denote by $I_{L_1}, I_{V_{1+}}, I_{A_1}, ..., I_{L_4}$ the parts of the integral \eqref{MajorSmallMinor}  over the corresponding closed line segment or  arc. Choose a function $\mu:\R_{\ge 0}\to \R_{\ge 0}$, $\lim_{t\to 0+} \mu(t)=0$, such that $|\phi(\xi)|\le \frac{\mu(t)}{|\xi|}$ for $|\xi|<t$, $\arg \xi\in B''$.  Then 
$$ I_{A_1} \le \frac{ length(A_1) \cdot \max_{\eta\in A_1} ( \phi(\eta) ) }{ r_\xi }  \le 2\pi\max_{\eta\in A_1} ( \phi(\eta) )  \ \le $$
$$ \le 2\pi\frac{\mu(2|\xi|)}{(|\xi|/2)} \ = \ o(\frac{1}{|\xi|}), $$
and similarly for $A_3$.

Now notice that on all other arcs and segments $|\xi-\eta|>|\xi|r_\xi \ge \frac{1}{3}|\xi| \sin\Delta\beta$. \\
Hence, e.g., 
$$ I_{L_2} \le \frac{ length(L_2)\cdot \max_{\eta\in L_2}{|\phi(\eta)|} }{\frac{1}{3}|\xi| \sin\Delta\beta } . $$ 
Since $length(L_2)< |\xi|$, and for $\eta\in L_2$ we have $|\phi(\eta)|\le \mu(|\xi|)/|\eta|$ and $|\eta|\ge |\xi|/3$, we get
$$ I_{L_2} \le \frac{3}{\sin\Delta\beta } \frac{\mu(|\xi|) }{ |\xi|/3} =  \frac{9}{\sin\Delta\beta } \frac{\mu(|\xi|) }{ |\xi|} = o(\frac{1}{|\xi|}), $$
and analogously for $V_{1\pm}$, $A_2$, $L_3$, $V_{3\pm}$. 

As for $I_{L_1}$ (and similarly for $I_{L_4}$), split $L_1$ into the union of two intervals $[\Re \xi  + r_\xi, \frac{3}{2}|\xi|]$ (which is of length $O(\xi)$) and $[\frac{3}{2}|\xi|; 2T]$. The integral over the former interval can be estimated as above, and let us show that
$$ |\xi|\int_{\frac{3}{2}|\xi|}^{2T}  \frac{\phi(\eta)}{\xi-\eta} d\eta \to 0 , \ \ \text{as} \ \xi\to 0.  $$
Analogously to \cite[Pr\'e II.5.2, p.166]{CNP}, change the variable $\eta=|\xi| t$, writing $\xi=|\xi|e^{i\theta}$, we can rewrite the statement as
$$ \int_{\frac{3}{2}}^{2T/|\xi|} \frac{|\xi| \phi(|\xi|t)}{(e^{i\theta}-t)} dt \to 0 , \ \ \text{as} \ \xi\to 0,  $$
which follows by the dominated convergence theorem. 



Adding together these estimates, we obtain the lemma. $\Box$

\begin{Lemma} \label{TrulySmall} Let $\varphi(h)$ be any representative of a small resurgent function in the direction $A$. Then, for any arc $A'\subset\subset A$ and any $\varepsilon>0$ there is a sectorial neighborhood $U$ of $0$ in the direction $A'$ such that $|\varphi(h)|<\varepsilon$, $\forall h \in U$. \ $\Box$ \end{Lemma}

\textsc{Proof.} After reducing the question to the case of only one singularity at the origin, use the proof of ~\cite[Pr\'e II.5.2, p.166]{CNP}. Ingredients of that proof have been detailed here, see lemma \ref{MajorSmallResurgFctn} and proposition \ref{Prop351}. $\Box$

\subsection{Minors} \label{SubsectionMinors}

Following  ~\cite[Pr\'e II.4]{CNP}, 
denote ${}^\flat {\pmb{\cal C}}(A)$ the subalgebra (with respect to convolution) of small mictofunctions in ${\pmb{\cal C}}(A)$. (It is indeed a subalgebra with respect to convolution. ) 

From now on consider microfunctions defined on a big arc $B$ . For an arc $A$ denote by ${\cal O}^0(A)$ the space sectorial germs of analytic functions at $0$ in the direction $A$. We can define a variation
$$ {\rm var} \ : \ {\pmb {\cal C}}(B) \to {\cal O}^0(\hat B)$$
as follows: take a microfunction $\phi$ at $0$, analytically continue it once counterclockwise around zero, obtain a microfunction $\tilde \phi$, and put ${\rm var} \ \phi := \phi - \tilde \phi$.

We will show that small microfunctions are specified by their variation (or its ``minor").

Denote by ${}^{\rm min} {\cal O}^{0}(\hat B)$ the space of germs of holomorphic functions in a sectorial neighborhood $V$ of direction $\hat B$ whose primitive tend to a finite limit when $\zeta\to 0$ in $B$. Denoting by $G(\zeta)$ this primitive and considering the family of continuous functions $f_r(e^{i\tau})=G(re^{i\tau})$ of $\tau$ indexed by $0<r<<1$, we see that for $\tau$ in a compact subinterval of $B$ the convergence for $r\to 0+$ can be made uniform. 

\begin{Prop} \label{Prop351} (\cite[Pr\'e II.4.2.1]{CNP}) The map $\var$ is an isomorphism ${}^\flat{\pmb{\cal C}}(B) \to {}^{\rm min}\O^0(\hat B)$. Denote by b\'emol its inverse $g\mapsto {}^\flat g$.
\end{Prop}

\textsc{Proof.} The proof we are going to present is ~\cite{CNP}'s proof changed and clarified using suggestions of this paper's two anonymous referees.

\textsc{Part I.}  Let us show that ${\rm var}$ send ${}^\flat{\pmb{\cal C}}(B)$ to ${}^{\rm min}\O^0(\hat B)$. Indeed, let $\varphi$ be represented by a function ${\bf \Phi}$ holomorphic in a sectorial neighborhood of $0$ in direction $B$ and satisfying $\Phi(\xi)=o(1/|\xi|)$ for $\xi\to 0$.

\begin{figure}[h]\includegraphics{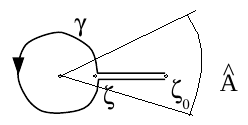} \caption{Contours in the proof of Proposition \ref{Prop351}} \label{ResDRWp80}
\end{figure}

Let $V$ be a sectorial neighborhood of $0$ in direction $\hat B$, let $\zeta_0\in V$, let $\gamma$ be a contour starting at $\zeta_0$ and going around $0$, as on the figure \ref{ResDRWp80}. 
  
We have 
$$\int_\gamma {\bf \Phi}(\xi)d\xi \ = \ \int_{C_\zeta} {\bf\Phi}(\xi) d\xi \ - \ \int_{[\zeta_0,\zeta]} ({\rm var \ }\varphi)(\xi) d\xi, $$
where $C_\zeta$ is a small circular contour starting from a point $\zeta$ close to $0$.

As $\int_\gamma {\bf \Phi}(\xi)d\xi$ is independent of $\zeta$ and the integral over $C_\zeta$ tends to zero when $\zeta\to 0$, the function $\zeta\mapsto \int_{C_\zeta} {\bf \Phi}(\xi)d\xi$ is a primitive of ${\rm var \ }\varphi$ that is finite at zero. 

\textsc{Part II.} Construction of the inverse map ``b\'emol".

Let $g\in {}^{\rm min} \O^0(\hat B)$. To construct a representative ${\pmb\Phi}$ of a microfunction ${}^\flat g$, let $G(\eta)$ be the primitive of $g$ that tends to $0$ when $\eta\to 0$ and choose $T>0$ so that $G(\xi)$ has no singularities for $0<|\xi|<3T$ and $\arg \xi\in B$. Let $\eta_1\in \C$, $\arg \eta_1\in \hat B$ and $|\eta_1|=2T$.
Put  
$$ {\pmb \Phi}(\zeta) \ = \  \frac{1}{2\pi i} \int_{[0,\eta_1]} \frac{G(\eta)}{(\eta-\zeta)^2}d\eta $$
for $\zeta$ away from the line segment $[0,\eta_1]$; for other values of $\zeta$ define the analytic continuation by the deformation of the integration path.
 
We need to show that ${\pmb\Phi}(\zeta)$ is $o(1/|\zeta|)$ for $\zeta\to 0$ uniformly in every angular sector compactly contained in $B$. 

For simplicity of notation let $\eta_1\in \R_{\ge 0}$. 
Similarly to the proof of lemma \ref{MajorSmallResurgFctn}, take subsectors $B'\subset\subset B''\subset\subset B$.  For $|\zeta|<T$, $\arg \zeta\in B'$  choose a path $\Gamma_\zeta$ on the complex plane of the variable $\eta$ it will consist from a path $\Gamma_\zeta^0$ from $\eta=0$ to $\eta=2|\zeta|$ contained in $\{ \eta: \ \arg\eta\in B''\}$ and staying a distance of order $|\zeta|$ from $\zeta$ (one can choose such a path similarly to $\Gamma_\xi$  from the proof of lemma \ref{MajorSmallResurgFctn}); and the line segment $[2|\zeta|,\eta_1]$. Similarly to the proof of lemma \ref{MajorSmallResurgFctn}, one shows that $\int_{\Gamma_\zeta^0}\frac{G(\eta)}{(\eta-\zeta)^2}d\eta$ is $o(1/|\zeta|)$ and it remains to show that
$$ \lim_{|\zeta|\to 0}\left[ \zeta \int_{2|\zeta|}^{\eta_1} \frac{G(\eta)}{(\eta-\zeta)^2} d\zeta \right] \ = \ 0. $$
As in ~\cite[Pr\'e II.5.2]{CNP}, changing the variable by $\eta=|\zeta| t$ transforms the question to showing
$$ \lim_{|\zeta|\to 0} \left[ \int_{2}^{\eta_1/|\zeta|} \frac{G(|\zeta|t)}{(t-2)^2} dt \right] = 0 $$
which follows by the dominated convergence theorem.

Finally, one checks using the Cauchy theorem that ${\rm var \ }[{\pmb\Phi}]=g$. $\Box$

It is easy to see from the above formulas that $\Phi$ is defined Riemann surface obtained from the Riemann surface of  $g$ by adding a branch point on its every sheet over $\eta_1$. \footnote{ \cite{CNP} say it is defined on the same Riemann surface as $g$.} The function $\Phi$ is called the \underline{adapted major} in \cite{CNP}, or more to the author's taste, the \underline{adapted representative} of our microfunction.

Suppose $[\phi],[\psi]$ are resurgent microfunctions, and their variation has no singularities on the line segment $[0,\eta]$.  Then ~\cite{CNP} write that for  $t\in(0,\eta]$ 
\begin{equation}  [\var (\phi *  \psi)](t)  \ = \ \int_{(0,t)} (\var \phi(\tau))(\var \psi(t-\tau)) d\tau  \label{CNPminorCon} \end{equation}
and for $t$ farther away from $0$ we might need to use analytic continuation. Following a suggestion of the referee, since we do not know that $\phi$ and $\psi$ are integrable at the origin, we will take \ref{CNPminorCon} to mean the following:
$$ [\var (\phi *  \psi)](t)  \ = \ \frac{d}{dt} \left( \int_{(0,t)} G(\tau) H(t-\tau) d\tau \right), $$ 
where $G$ and $H$ are primitives of $\var \phi$ and $\var \psi$, respectively.


\section{Substitution of a small resurgent function into a holomorphic function}

The goal of this section is to prove the following theorem:

\begin{Thm} \label{substituteThm} If $g(z_1,...,z_k)=\sum a_{j_1...j_k}z_1^{j_1}...z_k^{j_k}$ is a complex analytic function given around the origin by a convergent power series, and $\varphi_1(h),...,\varphi_k(h)$ are small resurgent functions, then the composition $g(\varphi_1,...,\varphi_k)$ is a resurgent function. 
\end{Thm}

In ~\cite{CNP} 
an analogous result is proven for $k=1$ and for the case of a resurgent function $\varphi_1(h)$ representable by a major with a single singularity at the origin; recall the difference between our and ~\cite{CNP}'s definition of a small resurgent function explained in section \ref{SmallResFunSec}. Our definition, in contrast to ~\cite{CNP}'s, includes, e.g., $\varphi(h)=h^2+e^{-1/h}$ as a small resurgent function; this function is not fully determined by the singularity of its major at the origin or by the minor corresponding to $h^2$; therefore ~\cite{CNP}'s method of proof by estimating iterated integrals of minors remains an important special case, but is no longer sufficient for the proof theorem \ref{substituteThm}.

A sketch of ~\cite{CNP}'s proof is presented in section \ref{kvarMinors}, together with the changes necessary to pass to the $k$ variable case. In section \ref{ReduceToMinor} we reduce our, more general, situation to the one treated by ~\cite{CNP}. 

We remark that theorem \ref{substituteThm} can be also derived by induction on the number of variables from a parameter-dependent version of proposition \ref{Prop9}.






\subsection{Generalization of ~\cite{CNP}'s construction of a composite function  to the case of $k$ variables.}
\label{kvarMinors}

First let us mention that ~\cite{CNP} work with a weaker definition of an endlessly continuable function. For them the function is endlessly continuable if it can be analytically continued along any path of length $L$ and angle variation $\delta$ avoiding a finite set $\Omega_{L,\delta}$. To our knowledge, resurgent functions in this weaker sense possess all useful properties of resurgent functions in the stronger sense. 

Following \cite[R\'es II.3]{CNP}, for a resurgent function consider the Riemann surface ${\cal S}$ of its major and the discrete filtered set $\Omega_*=\{ \Omega_{L} \}$ of its singularities, accessible along paths $\gamma$ of length $L$ with fixed starting point $\gamma(0)$.

Conversely, for any discrete filtered set $\Omega_*$, ~\cite{CNP} constructs a Riemann surface ${\cal S}(\Omega)$. 
We feel, however, that a clear exposition of the conditions that should be imposed on $\Omega_*$ in order for ${\cal S}(\Omega)$ to be an {\it endlessly continuable} Riemann surface, is still absent from the literature.

The sum of two discrete filtered set $\Omega'_*$ and $\Omega''_*$ denoted $(\Omega'+\Omega'')_*$ is defined as follows:\\
$(\Omega'+\Omega'')_{L}$ is the set of $\omega'+\omega''$ where $\omega'\in \Omega'_{L'}$, $\omega''\in \Omega''_{L''}$ and $L'+L''=L$. 

Let ${}^n\Omega_*$ denote the sum of $n$ copies of $\Omega_*$; let ${}^\infty\Omega_* : = \bigcup {}^n\Omega_*$ 

 If two resurgent functions have $\Omega'_*$ and $\Omega''_*$ as filtered discrete sets of their singularities, the singularities of their convolution are included in $(\Omega'+\Omega'')_*$, cf. ~\cite[R\'es II.3.1.5]{CNP}.

If ${}^\flat{\pmb f}$ is a small microfunction, denote by ${\pmb f}^{(-1)}$ the primitive of ${\pmb f}$ that vanishes at $0$.   
We will be using the concept of the convolution of minors from ~\cite[Pr\'e II.4.2.2]{CNP}.

Given small resurgent function ${}^\flat{\pmb f}_1,..., {}^\flat{\pmb f}_n$, the function  ${}^\flat{\pmb f}_1*...* {}^\flat{\pmb f}_n=:{}^\flat{\pmb g}$ is also small, and for $\zeta$ close to $0$ we have 
$$ {\pmb g}^{(-n-1)}(\zeta) = \int {\pmb f}_1^{(-1)}(s_1){\pmb f}_2^{(-1)}(s_2)...{\pmb f}_n^{(-1)}(s_n) ds_1 ds_2 ... ds_n, $$
where the integral is taken over the $n$-simplex defined by 
$$ \arg s_1 = ... = \arg s_n = \arg \zeta, \ \ \ \ |s_1|+|s_2|+...+|s_n|\le |\zeta|. $$

In order to obtain some estimates on the growth of ${\pmb g}^{(-n-1)}$ it will be shown that it is possible to define the continuation of ${\pmb g}^{(-n-1)}$ along any allowed path 
as an integral over an $n$-simplex obtained by a deformation of the initial simplex.

Let ${\cal S}$ be an endless Riemann surface and $\Omega_*$ the discrete filtered set of its singularities.  A sectorial neighborhood of $0$ in ${\cal S}$ is said to be {\it small} if all its points are $L$-accessible from $0$, with $L$  so small that $\Omega_{L}=\{ 0\}$;\\
a sectorial neighborhood of $0$ is {\it bounded} if it is the union of a small neighborhood of $0$ and a relatively compact subset of ${\cal S}$. 

If ${\pmb f}$ is a minor of a small resurgent function, then its primitive ${\pmb f}^{(-1)}$ is bounded on any bounded neighborhood of $0$.

Let ${\cal S}:={\cal S}(\Omega_*)$  be the Riemann surface on which ${\pmb f}_1,...,{\pmb f}_k$ are simultaneously defined. Let ${}^\infty{\cal S}:={\cal S}({}^\infty \Omega_*)$. 

The following corresponds to \cite{CNP}'s Key lemma 2. 

\begin{Lemma} \label{KeyLemma2}  There is an exhaustive family of bounded neighborhoods $(V_{\alpha,\infty} \subset {}^\infty{\cal S})$, a family of bounded neighborhoods $(V_\alpha\subset {\cal S})$, constants $C_\alpha >0$ , and for any minor ${\pmb f}$ of a small resurgent microfunction which is analytic on ${\cal S}$, there is a family of functions $\epsilon_\alpha$ defined on $\N$, s.th. $\epsilon_\alpha(n)>0$ and $\epsilon_\alpha(n)\to 0$ as $n\to +\infty$, so that if ${\pmb g}_{j_1...j_k}$ denotes the minor of $({{}^\flat\pmb f}_1)^{*j_1}*...*({{}^\flat\pmb f}_1)^{*j_k}$ with $j_1+...+j_k=n$, one has an estimate
\begin{equation} |{\pmb g}_{j_1...j_k}^{(-n-1)}|_{V_{\alpha,\infty}} \ \le \ \frac{1}{n!}  [  C_{\alpha}  ( \max_j |{\pmb f}^{(-1)}_j|_{V_{\alpha}})^{1/2} \epsilon_{\alpha} (n)]^n  \label{eq20} \end{equation}
\end{Lemma}

The proof is an easy modification of the one given in ~\cite{CNP}. Once we have established the counterpart of Key lemma, the rest of the proof goes along the same lines as in ~\cite{CNP}.

In particular, by the Cauchy integral formula, on a compact subset $K\subset V_{\alpha,\infty}$ that is separated by a distance $r_K$ from the branching poins of ${}^\infty{\cal S}$, we have an estimate
\begin{equation*} |{\pmb g}_{j_1...j_k}^{(-1)}|_{K} \ \le \  [  C_{\alpha}  r_K^{-1} ( \max_j |{\pmb f}^{(-1)}_j|_{V_{\alpha}})^{1/2} \epsilon_{\alpha} (n)]^n.   \end{equation*} 
As $\varepsilon_\alpha(n)\to 0$, the sequence of primitives ${\pmb g}_{j_1...j_k}^{(-1)}$ converges on compact sets faster than a geometric series with any positive ratio.  
Now, using ~\cite[R\'es 3.2.5]{CNP} or proposition \ref{Prop341}, we can find a sequence of majors $G_{j_1...j_k}$ corresponding to minors ${\pmb g}_{j_1...j_k}$ that also converges on compact sets of their Riemann surface faster than a geometric series with any positive ratio. Let us formulate this for the future use.

{\bf Lemma \ref{KeyLemma2}.A.} {\it Let $A$ be a small arc, $\alpha\in A$, $\varphi(h)$ a small resurgent function in the direction $A$ such that ${\bf s}_{\alpha+}^{-1} \bar{\cal L}\varphi$ consists a single microfunction at the origin. Then for the powers $[\varphi(h)]^n$, $n\ge 1$, there is a choice of majors $G_n$ that are bounded by geometric series with any a priori chosen ratio on compact subsets of their common Riemann surface.
}


\subsection{Reduction to the case of a small resurgent function with only one singularity of the major} \label{ReduceToMinor}

To simplify our notation, we will restrict ourselves to the case $k=1$ in this subsection.




Let $A$ be a small arc, $Sec(0,\check A)$ a closed sector with vertex $0$ in the direction $\check A$ containing the positive real direction; the isomorphism ${\bf s}_{0+}$ used in the proof below can just as well be consistently replaced with ${\bf s}_{0-}$. All convolutions of majors in this subsection will be taken with respect to a fixed contour $\Gamma$ adapted to $Sec(0,\check A)$; convolutions are then endlessly continued by deforming bounded segments of $\Gamma$ as explained in section \ref{ConvolutionSectn}.

Let $g(x)=\sum_{j=0}^\infty a_j x^j$ be a convergent series representing a holomorphic function $g(x)$ near the origin, and let $\varphi(h)$ be a small resurgent function represented by a major $F(\xi)+R(\xi)$ holomorphic in $Sec(0,\check A)$, where ${\bf s}_{0+} F(\xi)$ consists of a single microfunction at $\xi=0$  and  ${\bf s}_{0+}R(\xi)$ is supported inside the half-plane $\Re \xi>0$. Respectively, $\varphi(h)=\varphi_0(h)+r(h)$, where $\varphi_0={\cal L}F$, $r={\cal L}R$.

The construction we are going to describe has been greatly simplified using an idea of the anonymous referee. Notice that if the series $\sum_{j=0}^\infty a_j x^j$  has a nonzero radius of convergence, then so does $\sum_{j=0}^\infty C^j_k a_j x^{j-k}$ for any $k\ge 1$.   Let ${\pmb f}$ be the minor $\var F$, then by section \ref{kvarMinors} one obtains endlessly continuable minors ${\pmb h}_k=\sum C^j_k a_j [{\pmb f}^{*(j-k)}]$ for $k\ge 0$, and hence the corresponding  endlessly continuable majors $H_k={\bf s}_{0+}({}^\flat{\pmb h}_k)$. Then one can define a Mittag-Leffler sum of  majors $H_k * R^{*k}$ as in section \ref{MLS}, $\Psi(\xi)=\MLS_k H_k*R^{*k}$. 

This $\Psi(\xi)$ is our candidate for the major of $g(\varphi(h))$. Now  we need to take to calculate its Laplace transform of $\Psi$ and show that the result is equal to the  sum of the power series $\sum_j  a_j \varphi(h)^j$. Indeed, we have the following equalities  of resurgent functions:
$$  \begin{aligned}   {\cal L}\left[ \MLS_k H_k * R^{*k}  \right]  
&  \stackrel{(1)}{=} \ 
\sum_k {\cal L}\left[   H_k  *R^{*k}  \right]  \ = \ \cr 
&  = \ \sum_k ({\cal L} H_k) \cdot ({\cal L} R)^{k} \ = \ \cr 
& \stackrel{(2)}{=} \ \sum_{k=0}^\infty \left(\sum_{j=0}^\infty C_k^j a_j [\varphi_0(h)]^{j-k} \right) [r(h)]^k \ = \ \cr
& \stackrel{(3)}{=} \  \sum_{j,k} a_j C_j^k  [\varphi_0(h)]^{j-k} [r(h)]^k \ = \ \cr
& = \ \sum_j a_j (\varphi_0(h) + r(h))^j.  
\end{aligned}
$$
Here (1) holds by corollary \ref{CorLMLS}, (2) by ~\cite{CNP}'s construction recalled in section P\ref{kvarMinors}, and (3) by the Weierstrass' argument (cf, e.g., \cite[p.22]{Hi}) applicable since $\varphi_0(h)\to 0$ and $r(h)\to 0$   when $h\to 0$ in the given sector $A$ by lemma \ref{TrulySmall}; the double sum on the right-hand-side of (3) is an absolutely convergent double sum. 

This finishes the proof. The argument can be easily generalized to the case of $k$ variables. 
 



%


\subsection{Parameter-dependent version.}

The following is an easy parameter-dependent version of theorem \ref{substituteThm}.

\begin{Thm} Suppose $r(E,h)$ is a small resurgent function such that its major $R(E,\xi)$ is defined on one and the same Riemann surface and analytically dependent on $E$ for $E\in U$ a neighborhood if $0$. Suppose $f\in {\cal O}(D_\rho)$. Then we can choose a major $\Phi(E,\xi)$ of $f(r(E,h))$ for which the same is true. \end{Thm}

\textsc{Proof} is analogous to the one given for the parameter-independent case. The only modification is that in \ref{eq20} one has to replace $|{\pmb f}_j^{-1}|_{V_\alpha}$ by $|{\pmb f}_j^{-1}|_{V_\alpha\times \overline{D_{\rho'}}}$ for any $\rho'<\rho$. $\Box$

Here is an application of this theorem. A major of a solution of an $h$-differentail equation $P\psi=E_r h \psi$ can be chosen to analytically depend on $E_r$, and the same is therefore true for microfunctions -- formal solutions of the above equation. We can apply the above theorem to conclude that quotients of such microfunctions, in particular, formal monodromies can be represented by majors that holomorphically depend on $E_r$, and the same is then also true for the formal monodromy exponents. So, we have:

\begin{Cor} Formal monodromy exponents for an $h$-differential equation can be represented by majors holomorphically dependent on $E_r$. \end{Cor}


\section{Substitution of a small resurgent function for a holomorphic parameter of another resurgent function.}

Let us fix a small arc $A$ containing, for definiteness, the positive real direction. In this section by resurgent or small resurgent functions we mean (small) resurgent functions in the direction $A$.

\begin{Prop} \label{Prop9} Suppose $\varphi(E,h)$  is a resurgent function defined for every $E\in \Delta$, where $\Delta\subset \C$ is a disc of radius $r>0$ around the origin, and  satisfying the following property:\\
For all $E\in \Delta$ the majors $\Phi(E,\xi)$ are defined on the same Riemann surface ${\cal S}$ and $\Phi(E)$ is a holomorphic function on $E\times {\cal S}$. \\
Then for any small resurgent function $E(h)$, the function $\varphi(E(h),h)$ is resurgent. \end{Prop}

\textsc{Proof.} Construct the Riemann surface ${\cal T}$ for the major of $\varphi(E(h),h)$ as ${\cal S}*{}^{\infty}{\cal S}_E$ where ${\cal S}_E$ is the Riemann surface of ${\tilde E}(\xi)$

Let us construct the major for $\varphi(E(h),h)$. Let ${\tilde E}(\xi)=G + R$ where $G$ is repesented by only one small resurgent microfunction ${}^\flat {\pmb g}$ and $R$ has all its first-sheet singularities to the right of the imaginary axis. 

Without loss of generality, assume that $\Phi(E,\xi)$ is holomorphic in the sector $Sec(0,\check A)$.

a) By the lemma \ref{KeyLemma2}.A, we can choose representatives $G_n \in G^{*n} \mod {\cal O}(\C)$ that go to zero faster  than any geometric series (here convolutions are taken along any path adapted to $Sec(0,\check A)$). Therefore, by Proposition \ref{Prop331} we can take representatives $G_n \in G^{*n} \mod {\cal O}(\C)$ and a contour $\Gamma$ consisting of two half-lineas and adapted to $Sec(0,\check A)$, so that $G_n$ are bounded ``outside of $\Gamma$" by a geometric series, i.e. $|G_n|\le C q^{n}$ in that domain, and such that $|G_n|\le M_K q^{n}$ on any compact set $K\subset{\cal T}$, where $q=\frac{1}{2}\min\{1,r\}$. 

The path $\Gamma$ splits the Riemann surface ${\cal T}$  into two closed parts: $U_{-}$ which is the first sheet of ${\cal T}$ outside of $\Gamma$, i.e. containing a sectorial neighborhood of infinity in the direction $\check A$, and the rest, denoted $U_+$, so that $U_{-}\cap U_{+}=\Gamma$.  For a compact set $K\subset {\cal T}$ denote $K\cap U_{\pm} = K_{\pm}$.

b) By using a parameter-dependent version of lemma \ref{Lemma21}, after possibly shrinking $\Delta$, we can assume that $|\Phi(E,\xi)|<e^{-|\xi|}$ along $\Gamma$. Then $\Phi_n(\xi)=\frac{1}{n!}\frac{\partial^n \Phi(E,\xi)}{\partial E^n}$  will satisfy $\left| \Phi_n \right| < C Q^n e^{-|\xi|}$ on $\Gamma$ for some $0<Q<1$.
 Further, for compact subsets $K\subset {\cal S}$, $\sup_{\xi\in K_+\cup U_-} |\Phi_n| < C_K r^{-n}$.

c) To show that $\sum_{n=0}^\infty C_n^j \int_\Gamma G_n(\xi-\eta)\Phi_n(\eta)d\eta$ converges uniformly on compact subsets faster than a geometric series, one imitates the proof of Proposition ~\ref{Prop331}. 

%

d) Now we need to show that this major is really the major of the composite function, i.e. we need to verify the equality 
$$ {\cal L} \left[\MLS_{j=0}^\infty R^{*j} * \left(\sum_{n=0}^\infty C^j_n \int_{\Gamma} G_n(\xi-\eta) * \Phi_n(\eta) d\eta 
\right) \right] \ = \  \sum \frac{1}{n!} \frac{\partial^n \phi}{\partial E^n} [E(h)]^n. $$

In the first sum $\MLS$ and ${\cal L}$ can be interchanged and we get
$$ \sum_{j=0}^\infty  {\cal L}[R]^{j} \cdot {\cal L} \left[\sum_{n=0}^\infty C^j_n \int_{\Gamma} G_n(\xi-\eta) * \Phi_n(\eta) d\eta  \right].$$

Since  the majors inside the sum converge uniformly and faster than geometric series with some fixed ratio on compact sets, ${\cal L}$ and the infinite sum can be interchanged, after which the proposition is proven. $\Box$


\section{Application to quantum resurgence.} \label{ApplicnToExistence}

Consider the Schr\"odinger operator 
$$ P \ = \ -h^2\partial_q^2 + V(q,h), $$
where $V(q,h)$ is analytic in $q$ on the whole complex plane and polynomial in $h$. Denote by $\tilde P$ the Laplace-transformed operator, i.e,
$$ \tilde P \ = \ -\partial_\xi^{-2}\partial_q^2 + V(q,\partial_\xi^{-1}), $$
and by $\hat h$ is a major of $h$, e.g. $\hat h(\xi)=\frac{1}{2\pi i} \log \xi$.

Suppose we know that $\psi$ is a resurgent solution  of the differential equation
$$ P\psi(E_1,q) = hE_1\psi(E_1,q) \ \ \ \mod {\cal E}^{-\infty}$$ 
for every $E_1\in \C$.  
Here ${\cal E}^{-\infty}$ stands for functions of $h$ that are $<C_{a,K}e^{-a/h}$ for small $|h|$ when other parameters range over a compact set. 
We will assume that the majors $\Psi(E_1,\xi,q)$ are holomorphic with respect to $E_1$, are defined on the same Riemann surface and satisfy
\begin{equation} {\tilde P}\Psi(E_1,q,\xi) = {\hat h}*E_1\Psi(E_1,q,\xi) \ \ \ \mod{{\cal O}}(U_1\times U_2\times \C), \label{whatSolusSatisfy} \end{equation}
where $U_1$ is an open neighborhood of $0$ in $\C$ and $q\in U_2$, an open subset of the universal cover $\tilde \C$ of $\C\backslash\{\text{turning pts}\}$ and we will assume $U_2$ to be relatively compact. (Turning points are those $q$ for which $V(q,h)=O(h)$.) Such are the properties we expect of solutions produced, say, by Shatalov-Sternin method, ~\cite{ShSt}.

Differentiating both sides of \eqref{whatSolusSatisfy} with respect to $E_1$ at $E_1=0$, obtain
$$ {\tilde P} \left.\frac{\partial^n\Psi}{\partial E_1^n}\right|_{E_1=0} \ = \ {\hat h}*\left.n\frac{\partial^{n-1}\psi}{\partial E_1^{n-1}}\right|_{E_1=0} \ \ \ \mod {{\cal O}}(U_1\times U_2\times \C).$$

As $\Psi(E_1,q,\xi)=\sum E_1^n \left.\frac{1}{n!}\frac{\partial^n \Psi}{\partial E_1^n}\right|_{E_1=0}$ converges uniformly on compact sets of the Riemann surface of $\Psi$ and is holomorphic with respect to $E_1$, we see that $\frac{1}{n!}\frac{\partial^n \Psi}{\partial E_1^n}$ can be estimated on compact sets of that Riemann surface by geometric series with the same ratio. Hence, by the previous section, we may substitute a small resurgent function $E(h)$ for $E_1$ and obtain a resurgent function.

We know that  $P\psi(E_1,h) = hE_1\psi(E_1,h) \mod {\cal E}^{-\infty}$ for $E_1\in \C$. If they were equal as functions, not as classes modulo ${\cal E}^{-\infty}$, then there would be no issue substituting $E(h)$ into this equality. As it is, we need an additional argument in order to show: 
$$ P\psi(E(h),q) = hE(h) \psi(E(h),q) \mod {\cal E}^{-\infty}.$$

On the level of majors it boils down to showing 
\begin{equation} {\tilde P} \sum {\tilde E}^{*n} * \frac{1}{n!} \left.\frac{\partial^n \Psi}{\partial E_1^n}\right|_{E_1=0} \ = \ {\hat h}*{\tilde E} * \sum {\tilde E}^{*n} * \frac{1}{n!} \left.\frac{\partial^n \Psi}{\partial E_1^n}\right|_{E_1=0} \mod {{\cal O}}(U_2\times \C). \label{eq22} \end{equation}

We can interchange ${\tilde P}$ and the infinite sum on the left. Indeed, for the convolution with ${\hat h}^2$ and $V(q,\hat h)$ this follows by proposition \ref{Prop43rev}. Infinite sum and $\frac{\partial}{\partial q}$ can be interchanged because on compacts with respect to $q$ the terms of the series can be assumed  bounded by $Q^j e^{-|\xi|}$, for some $0<Q<1$, along the infinite branches of $\Gamma$, and since everything is analytic in $q$, so are the $q$-derivatives of majors.

 We can interchange ${\tilde E}*$ and the infinite sum on the right. Indeed, it is shown in the proof of proposition \ref{Prop9} that majors of ${\tilde E}^{*n}*\frac{1}{n!}\frac{\partial^n\Psi}{\partial E_1^n}$ have representative converging faster than a geometric series on compact subsets. Therefore, by Prop. \ref{Prop43rev}, we can interchange ${\tilde E}*$ and the infinite sum.

So, to show \eqref{eq22}, since interchanging $\tilde P$ and ${\tilde E}*$ with the infinite sums in question is legal,  it is enough to show that 
$$ \sum {\tilde E}^{*n} * \frac{1}{n!}{\tilde P}\frac{\partial^n\Psi}{\partial E_1^n} \ = \ \sum {\hat h}*{\tilde E}^{*(n+1)}*\frac{1}{n!}\frac{\partial^n\Psi}{\partial E_1^n} \mod {{\cal O}}(U\times \C)$$
which follows because 
$$ {\tilde P}\frac{\partial^n\Psi}{\partial E_1^n} =  n {\hat h}*\frac{\partial^{n-1}\Psi}{\partial E_1^{n-1}} \mod {{\cal O}}(U\times \C) . $$
$\Box$


\vspace{2cm}

{\bf Acknowledgements} 
The author would like to thank his advisor Boris Tsygan for a wonderful graduate experience and his dissertation committee members Dmitry Tamarkin and Jared Wunsch for their constant support in his study and research. 
 We greatly appreciate thoughtful remarks of this paper's two anonymous referees.
Valuable suggestions were also made by M.Aldi, J.E.Andersen, K.Burns, K.Costello, E.Delabaere, S.Garoufalidis, E.Getzler, B.Helffer, A.Karabegov, S.Koshkin, Yu.I.Manin, G.Masbaum, D.Nadler, W.Richter, K.Vilonen, F.Wang, E.Zaslow, and M.Zworski.
This work was partially supported by the NSF grant DMS-0306624 and the WCAS Dissertation and Research Fellowship.

\end{document}